\documentclass[11pt]{article}
\usepackage{latexsym,amssymb}
\usepackage{amsmath}
\usepackage{color}

\newtheorem{Theorem}{\bf Theorem}[section]
\newtheorem{Lemma}{\bf Lemma}[section]
\newtheorem{Proposition}{\bf Proposition}[section]
\newtheorem{Corollary}{\bf Corollary}[section]
\newtheorem{Remark}{\bf Remark}[section]
\newtheorem{Example}{\bf Example}[section]
\newtheorem{Definition}{\bf Definition}[section]

\newenvironment{theorem}{\begin{Theorem}$\!\!\!$}{\end{Theorem}}
\newenvironment{lemma}{\begin{Lemma}$\!\!\!$}{\end{Lemma}}
\newenvironment{proposition}{\begin{Proposition}$\!\!\!$}{\end{Proposition}}

\newenvironment{remark}{\begin{Remark}$\!\!\!$}{\end{Remark}}

\newenvironment{definition}{\begin{Definition}$\!\!\!$}{\end{Definition}}

\numberwithin{equation}{section}

\begin{document}

\title{A supercritical scalar field equation with a forcing term}
\author{
\qquad\\
Kazuhiro Ishige, Shinya Okabe and Tokushi Sato
}
\date{}
\maketitle
\begin{abstract}
This paper is concerned with the elliptic problem for a scalar field equation with a forcing term
\begin{equation}
\tag{P}-\Delta u+u=u^p+\kappa\mu \quad\mbox{in}\quad{\bf R}^N,
\quad
u>0\quad\mbox{in}\quad{\bf R}^N,
\quad
u(x)\to 0\quad\mbox{as}\quad |x|\to\infty, 
\end{equation}
where $N\ge 2$, $p>1$, $\kappa>0$ and $\mu$ is a Radon measure in ${\bf R}^N$ with a compact support. 
Under a suitable integrability condition on $\mu$, 
we give a complete classification of the solvability of problem~(P) with $1<p<p_{JL}$. 
Here $p_{JL}$ is the Joseph-Lundgren exponent defined by 
$$
 p_{JL} :=\infty\quad\mbox{if}\quad N\le 10,
 \qquad
p_{JL}:=\frac{(N-2)^2-4N+8\sqrt{N-1}}{(N-2)(N-10)}\quad \text{if} \quad N\ge 11.
$$
\end{abstract}
\vspace{25pt}
\noindent Addresses:

\smallskip
\noindent K.~I.:  Graduate School of Mathematical Sciences, The University of Tokyo,
3-8-1 Komaba, Meguro-ku, Tokyo 153-8914, Japan.\\
\noindent 
E-mail: {\tt ishige@ms.u-tokyo.ac.jp}\\

\smallskip
\noindent S. O.:  Mathematical Institute, Tohoku University,
Aoba, Sendai 980-8578, Japan.\\
\noindent 
E-mail: {\tt okabes@m.tohoku.ac.jp}\\

\smallskip
\noindent 
T. S.: Miyagi University of Education,
Aoba, Sendai 980-0845, Japan.\\
\noindent 
E-mail: {\tt tokusi-s@staff.miyakyo-u.ac.jp}\\
\vspace{20pt}
\vspace{3pt}
\newpage
%%%%%%%%%%%%%%%%%%%%%%%%%%%%%%%%%%%%%
%%%%%%%%%%%%%%%%%%%%%%%%%%%%%%%%%%%%%
\section{Introduction}
%%%%%%%%%%%%%%%%%%%%%%%%%%%%%%%%%%%%%
%%%%%%%%%%%%%%%%%%%%%%%%%%%%%%%%%%%%%
This paper is concerned with the solvability 
of the elliptic problem for a scalar field equation with a forcing term
\begin{equation}
\label{eq:1.1}
\left\{
\begin{array}{ll}
-\Delta u+u=u^p+\kappa\mu & \quad\mbox{in}\quad{\bf R}^N,\vspace{3pt}\\
u>0 & \quad\mbox{in}\quad{\bf R}^N,\vspace{3pt}\\
u(x)\to 0 & \quad\mbox{as}\quad |x|\to\infty,
\end{array}
\right.
\end{equation}
where $N\ge 2$, $p>1$, $\kappa>0$ and 
$\mu$ is a nontrivial (nonnegative) Radon measure in ${\bf R}^N$ with a compact support. 
In particular, we are interested in problem~\eqref{eq:1.1} in the supercritical case $p>p_S$, 
where
$$
p_S:=\infty\quad\mbox{if}\quad N=2,\qquad p_S:=\frac{N+2}{N-2}\quad\mbox{if}\quad N\ge 3. 
$$
In general, 
the existence of the solutions of elliptic problems with supercritical nonlinearity is widely open 
since it is difficult to find the Sobolev embedding fitting suitably to a weak formulation of the solutions 
and many direct tools of calculus of variations are not applicable. 
See \cite{delPino01} and \cite{delPino02}, which include a nice survey and recent progresses 
for supercritical elliptic problems. 

In this paper, under a suitable integrability condition on $\mu$, 
we prove the existence of the critical constant $\kappa^*>0$ in the following sense (see Theorem~\ref{Theorem:1.1}): 
\begin{itemize}
  \item[(a)] 
  Problem~\eqref{eq:1.1} possesses a solution if $0<\kappa<\kappa^*$;
  \item[(b)] 
  Problem~\eqref{eq:1.1} possesses no solutions if $\kappa>\kappa^*$. 
\end{itemize}
Furthermore, we show the following result (see Theorem~\ref{Theorem:1.2}), 
which is the main result of this paper: 
\begin{itemize}
  \item[(c)] Let $1<p<p_{JL}$ and $\kappa=\kappa^*$. 
  Then problem~\eqref{eq:1.1} possesses a unique solution. 
\end{itemize}
Here 
$$
p_{JL} :=\infty\quad\mbox{if}\quad N\le 10,
\qquad
p_{JL}:=\frac{(N-2)^2-4N+8\sqrt{N-1}}{(N-2)(N-10)}\quad \text{if} \quad N\ge 11. 
$$
The exponent $p_{JL}$ is called the Joseph-Lundgren exponent (see \cite{JL}) 
and $p_{JL}>p_S$ for $N\ge 3$. 
It is a well-known critical exponent appearing in the study of the bifurcation structure 
of the radially symmetric solutions and the stability of the solutions of Lane-Emden-Fowler equation $\Delta u+u^p=0$  
(see e.g., \cite{Farina}, \cite{GNW01}, \cite{GNW02}, \cite{JL}, \cite{Miyamoto}, \cite{PY} and references therein). 
Throughout the proof of assertion~(c) 
we give a new characterization of the Joseph-Lundgren exponent $p_{JL}$. 
\vspace{5pt}

We recall some results closely related to this paper. 
Deng and Li \cite{DL01, DL02} proved assertions~(a) and (b) in $H^1({\bf R}^N)$ 
under the assumption that 
$$
\mu\in H^{-1}({\bf R}^N)\quad\mbox{and}\quad |x|^{N-2}\mu\in L^\infty({\bf R}^N).
$$ 
Furthermore, they proved assertion~(c) in the case of $1<p\le p_S$ and the following: 
\begin{itemize}
  \item[(d)]
  Assume either $1<p<p_S$ or $p=p_S$ with $3\le N\le 5$. 
  Then problem~\eqref{eq:1.1} possesses at least two solutions in $H^1({\bf R}^N)$ if $0<\kappa<\kappa^*$;
  \item[(e)]
  Let $p=p_S$ and $N\ge 6$. Under a suitable symmetric condition on $\mu$, 
  problem~\eqref{eq:1.1} possesses a unique solution in $H^1({\bf R}^N)$ 
  for all sufficiently small $\kappa>0$. 
\end{itemize}
On the other hand, 
the third author of this paper and Naito \cite{NS01} 
considered problem~\eqref{eq:1.1} with 
\begin{equation}
\label{eq:1.2}
\mu=\sum_{j=1}^n c_j\delta_{a_j},
\end{equation}
where $n\in\{1,2,\dots\}$, $c_j>0$, $a_j\in{\bf R}^N$ and $\delta_{a_j}$ is the Dirac measure supported at $a_j$. 
They proved assertions (a), (b), (c) and (d) in the case of  
\begin{equation}
\label{eq:1.3}
1<p<\infty\quad\mbox{if}\quad N=2,
\qquad
1<p<\frac{N}{N-2}\quad\mbox{if}\quad N\ge 3.
\end{equation}
We remark that problem~\eqref{eq:1.1} with \eqref{eq:1.2} possesses no solutions if $p\ge N/(N-2)$. 
There are many related results on assertions~(a)--(e). 
See e.g., \cite{BWZ}, \cite{CZ01}--\cite{Figu}, \cite{DL01}--\cite{DPW}, \cite{Hirano}--\cite{J},  
\cite{NS01}, \cite{NS02}, \cite{SW}, \cite{W}--\cite{ZZ} and references therein. 
However, unfortunately, they are not applicable to the proof of assertion~(c) in the supercritical case 
even if $\mu\in C({\bf R}^N)$ and $\mu$ has a compact support. 

In this paper we prove assertions~(a)--(c) 
and give a complete classification of the solvability of problem~\eqref{eq:1.1} 
in the case of $1<p<p_{JL}$. 
As far as we know, 
there are no available results for 
complete classifications of the solvability of problem~\eqref{eq:1.1} in the supercritical case. 
Multiple existence of solutions concerning assertions~(d) and (c) will be discussed in a forthcoming paper. 
\vspace{5pt}

We introduce some notation and formulate a definition of solutions of \eqref{eq:1.1}. 
For $x\in{\bf R}^N$ and $r>0$, let $B(x,r):=\{y\in{\bf R}^N\,:\,|x-y|<r\}$. 
Define
\begin{equation*}
\begin{split}
 & C_0({\bf R}^N):=\Big\{f\in C({\bf R}^N)\,:\,\lim_{|x|\to\infty} f(x)=0\Big\},\\
 & L^r_{\rm c}({\bf R}^N):=\{f\in L^r({\bf R}^N)\,:\,\mbox{$f$ has a compact support in ${\bf R}^N$}\},
\end{split}
\end{equation*}
where $1\le r\le\infty$. 
We denote by $G$ the fundamental solution of $-\Delta v+v=0$  in ${\bf R}^N$, that is
\begin{equation}
\label{eq:1.4}
G(x):=\frac{1}{(2\pi)^{N/2}|x|^{(N-2)/2}}K_{(N-2)/2}(|x|),
\qquad x\in{\bf R}^N\setminus\{0\},
\end{equation}
where $K_{(N-2)/2}$ is the modified Bessel function of order $(N-2)/2$. 
%%%%%%%%%%%%%%%%%%%%%%%%%%%%%%%%%%%%%%%%%%%%%%%
\begin{definition}
\label{Definition:1.1} 
Let $\mu$ be a Radon measure in ${\bf R}^N$, $\kappa>0$ and $1<p\le q<\infty$. 
\begin{enumerate}
\item[{\rm (i)}] 
We say that $u$ is a $(C_0+L^q_{\rm c})$-solution of \eqref{eq:1.1} if 
$u\in C_0({\bf R}^N)+L^q_{\rm c}({\bf R}^N)$ and 
$u$ satisfies 
$$
u(x)=[G*u^p](x)+\kappa[G*\mu](x)>0 \quad \text{for almost all} \quad x\in{\bf R}^N. 
$$
\item[{\rm (ii)}]
We say that $u$ is a $(C_0+L^q_{\rm c})$-supersolution of \eqref{eq:1.1} if 
$u\in C_0({\bf R}^N)+L^q_{\rm c}({\bf R}^N)$ and 
$u$ satisfies 
$$
u(x)\ge[G*u^p](x)+\kappa[G*\mu](x)>0 \quad \text{for almost all} \quad x\in{\bf R}^N. 
$$
\item[{\rm (iii)}]
We say that $u$ is a minimal $(C_0+L^q_{\rm c})$-solution of \eqref{eq:1.1} 
if $u(x)\le v(x)$ for almost all $x\in{\bf R}^N$ for any $(C_0+L^q_{\rm c})$-solution $v$ of \eqref{eq:1.1}. 
\end{enumerate}
\end{definition}

Now we are ready to state our results of this paper. 
Theorem~\ref{Theorem:1.1} is concerned with assertions~(a) and (b). 
\begin{theorem}
\label{Theorem:1.1} 
Let $\mu$ be a nontrivial Radon measure in ${\bf R}^N$ with $\mbox{supp}\,\mu\subset B(0,R)$ for some $R>0$. 
Let $p>1$ and assume that 
\begin{equation}
\label{eq:1.5}
G*\mu\in L^q({\bf R}^N)
\,\,\,\,\mbox{for some 
$\displaystyle{q>\max\biggr\{p,\frac{N(p-1)}{2}\biggr\}}$}.
\end{equation}
Then there exists $\kappa^{*}\in(0,\infty)$ with the following properties$\colon$  
\begin{itemize}
  \item[{\rm (i)}]
  If $1<\kappa<\kappa^{*}$, then problem~\eqref{eq:1.1}  
  possesses a minimal $(C_0+L^q_{\rm c})$-solution $u^\kappa$. 
  Furthermore, $u^\kappa(x)=O(G(x))$ as $|x|\to\infty${\rm ;} 
  \item[{\rm (ii)}]
  If $\kappa>\kappa^{*}$, then problem~\eqref{eq:1.1} possesses no $(C_0+L^q_{\rm c})$-solutions. 
\end{itemize}
\end{theorem}
\begin{remark}
\label{Remark:1.1}
We give some comments on assumption~\eqref{eq:1.5}.
\newline
{\rm (i)} Let $\mu\in H^{-1}({\bf R}^N)$. Then $G*\mu\in H^1({\bf R}^N)$, which together with the Sobolev embedding 
implies that $\mu$ satisfies \eqref{eq:1.5} in the case of $1<p<p_S$. 
\vspace{3pt}\newline
{\rm (ii)} Let $\mu$ satisfy \eqref{eq:1.2}. Then condition~\eqref{eq:1.5} holds 
under assumption~\eqref{eq:1.3}. 
\end{remark}
Due to Remark~\ref{Remark:1.1}, 
Theorem~\ref{Theorem:1.1} is somewhat new even in the subcritical case 
(compare with \cite{BWZ, CP, DL01, DL02, DPW, NS01}). 
The proof of Theorem~\ref{Theorem:1.1} is based on 
the construction of approximate solutions and the supersolution-subsolution method.  

In Theorem~\ref{Theorem:1.2}  
we show the unique solvability of problem~\eqref{eq:1.1} with $\kappa=\kappa^*$ in the case of $1<p<p_{JL}$. 
\begin{theorem}
\label{Theorem:1.2} 
Let $1<p<p_{JL}$ and 
assume the same conditions as in Theorem~{\rm\ref{Theorem:1.1}}. 
Then problem~\eqref{eq:1.1} with $\kappa=\kappa^*$ possesses a unique $(C_0+L^q_{\rm c})$-solution. 
\end{theorem}
The main ingredient in the proof of the existence of the solution with $\kappa=\kappa^*$ 
is to prove uniform local estimates of $\{w^\kappa\}_{0<\kappa<\kappa^*}$, 
where $w^\kappa:=u^\kappa-U_{j_*}^\kappa\ge 0$, 
$u^\kappa$ is the minimal solution given in Theorem~\ref{Theorem:1.1}  
and $U_{j_*}^\kappa$ is an approximate solution to \eqref{eq:1.1} (see \eqref{eq:3.2} and \eqref{eq:3.7}). 
Here $w^\kappa\in H^1({\bf R}^N)$ and it is a weak solution of a nonlinear elliptic problem (see \eqref{eq:4.3}). 
Applying elliptic regularity theorems to $\{(w^\kappa)^{1/2\nu}\}_{0<\kappa<\kappa^*}$, where $\nu\in(0,1)$,
with the aid of a delicate inequality (see Lemma~\ref{Lemma:5.1}),  
we obtain a uniform local estimate of $\{(w^\kappa)^{1/2\nu}\}_{0<\kappa<\kappa^*}$, 
instead of $\{w^\kappa\}_{0<\kappa<\kappa^*}$. 
This argument gives a new characterization of $p_{JL}$. 
Indeed, the argument requires to find $\nu\in(0,1)$ satisfying 
\begin{equation}
\label{eq:1.6}
4\nu(1-\nu)p>1,
\qquad
\frac{p_S+1}{2\nu}>\frac{N}{2}(p-1).
\end{equation}
The existence of $\nu\in(0,1)$ satisfying \eqref{eq:1.6} is equivalent to $1<p<p_{JL}$ (see Lemmas~\ref{Lemma:5.4} and \ref{Lemma:5.5}). 
Consequently, in the case of $1<p<p_{JL}$, 
we obtain a uniform local estimate of $\{(w^\kappa)^{1/2\nu}\}_{0<\kappa<\kappa^*}$ for some $\nu\in(0,1)$. 
Furthermore, we apply elliptic regularity theorems again to prove that $u^\kappa$ converges, as $\kappa\to\kappa^*$, 
to a solution~$u^{\kappa^*}$ of \eqref{eq:1.1} with $\kappa=\kappa^*$. 
The proof of the uniqueness of the solution with $\kappa=\kappa^*$ is by contradiction 
and the construction of supersolutions. 
\vspace{3pt}

The rest of this paper is organized as follows. 
In Section~2 we prepare some inequalities of the fundamental solution~$G$ and recall two lemmas on eigenvalue problems. 
In Section~3 we construct approximate solutions to \eqref{eq:1.1} and obtain some estimates of 
the approximate solutions.  
Section~4 is devoted to the proof of Theorem \ref{Theorem:1.1}. 
In Sections~5 and 6 we prove Theorem~\ref{Theorem:1.2}. 
%%%%%%%%%%%%%%%%%%%%%%%%%%%%%%%%%%%%%
%%%%%%%%%%%%%%%%%%%%%%%%%%%%%%%%%%%%%
\section{Preliminaries} \label{section:2}
%%%%%%%%%%%%%%%%%%%%%%%%%%%%%%%%%%%%%
%%%%%%%%%%%%%%%%%%%%%%%%%%%%%%%%%%%%%
In this section we recall some properties on the fundamental solution $G=G(x)$. 
In what follows, 
for any nonnegative functions $f$ and $g$ in ${\bf R}^N$, we say that 
$f(x)\asymp g(x)$ as $x\to 0$
if there exists $c>0$ such that $c^{-1}g(x)\le f(x)\le cg(x)$ in a neighborhood of $0$.  
Similarly, we say that  
$f(x)\asymp g(x)$ as $x\to\infty$
if there exists $c>0$ such that $c^{-1}g(x)\le f(x)\le cg(x)$ in a neighborhood of the space infinity.  
%%%%%%%%%%%%%%%%
\subsection{Fundamental solution~$G$}
%%%%%%%%%%%%%%%%
We collect some properties of the fundamental solution $G$. 
It follows from \eqref{eq:1.4} that 
\begin{equation}
\label{eq:2.1}
\begin{split}
 & G(x)\asymp 
\left\{
\begin{array}{ll}
|x|^{-(N-2)} & \mbox{if}\quad N\ge 3,\vspace{3pt}\\
-\log|x| & \mbox{if}\quad N=2,
\end{array}
\right.
\quad\mbox{as}\quad |x|\to 0,\\
 & G(x)\asymp |x|^{-\frac{N-1}{2}}e^{-|x|}
\qquad\qquad\qquad\quad\,\mbox{as}\quad |x|\to\infty. 
\end{split}
\end{equation}
By the H\"older inequality, 
the Hardy-Littlewood-Sobolev inequality and the Sobolev inequality 
we have the following properties (see \cite[Appendix]{NS01}): 
\begin{itemize}
  \item[(G1)] 
  For $r\in[1,N/(N-2))$, there exists $C_r>0$ such that 
  $$
  \|G*v\|_{L^r({\bf R}^N)}\le C_r\|v\|_{L^1({\bf R}^N)},
  \qquad
  v\in L^1({\bf R}^N);
  $$ 
  \item[(G2)]  
  For $r\in(1,N/2)$, there exists $C_r'>0$ such that 
  $$
  \|G*v\|_{L^{r'}({\bf R}^N)}\le C_r'\|v\|_{L^r({\bf R}^N)},
  \qquad
  v\in L^r({\bf R}^N),
  $$
  where $1/r'=1/r-2/N$;
  \item[(G3)] 
  Let $r>N/2$. Then 
  $$
  G*v\in C_0({\bf R}^N)\quad\mbox{and}\quad \|G*v\|_{L^\infty({\bf R}^N)}\le C\|v\|_{L^r({\bf R}^N)}
  $$ 
  for $v\in L^r({\bf R}^N)$;
  \item[(G4)]  
  For $v\in L^1({\bf R}^N)\cap L^r({\bf R}^N)$ with some $r>N/2$, then 
  $$
  G*v\in C_0({\bf R}^N)\cap L^1({\bf R}^N)\cap H^1({\bf R}^N).
  $$
\end{itemize}
Let 
\begin{align*}
g := G * \chi_{B(0,1)}, 
\end{align*}
where $\chi_{B(0,1)}$ denotes the characteristic function of the ball $B(0,1)$. 
Then  
\begin{equation}
\label{eq:2.2}
g\in C^\infty,
\quad
g(x)>0\quad\mbox{in}\quad{\bf R}^N,
\quad
g(x)\asymp G(x)\quad\mbox{as}\quad |x|\to\infty,
\quad
\frac{|\nabla g|}{g}\in BC({\bf R}^N).
\end{equation}
Furthermore, for any $\sigma>1$, 
there exists a constant $C>0$ such that 
\begin{equation}
\label{eq:2.3}
0<[G*g^\sigma](x)\le Cg(x)\quad\mbox{for}\quad x\in{\bf R}^N. 
\end{equation}
%%%%%%%%%%%%%%%%
\subsection{Eigenvalue problem}
%%%%%%%%%%%%%%%%
We recall two lemmas of the eigenvalue problem
\begin{equation}
\label{eq:2.4}
-\Delta\phi+\phi=\lambda a(x)\phi\quad\mbox{in}\quad{\bf R}^N,
\qquad
\phi\in H^1({\bf R}^N),
\end{equation}
where $a\in L^{N/2}({\bf R}^N)\cap L^r({\bf R}^N)$ for some $r>N/2$ 
and $a(x)>0$ for almost all $x\in{\bf R}^N$. 
See \cite[Lemmas~B2 and B3]{NS01}. 
\begin{lemma}
\label{Lemma:2.1}
Then problem~\eqref{eq:2.4} has the first eigenvalue $\lambda_1>0$ 
and the corresponding eigenfunction~$\phi_1$ with $\phi_1>0$ in ${\bf R}^N$. 
Furthermore, 
$$
\lambda_1=\inf\left\{
\|\psi\|_{H^1({\bf R}^N)}^2\,\biggr/\,\int_{{\bf R}^N}a\psi^2\,dx\,:\,
\psi\in H^1({\bf R}^N),\,\,\int_{{\bf R}^N}a\psi^2\,dx\not=0\right\}.
$$
\end{lemma}
\begin{lemma}
\label{Lemma:2.2}
Let $\lambda_1$ be the first eigenvalue to problem~\eqref{eq:2.4} 
and assume that $\lambda_1>1$. 
Then, for any $f\in H^{-1}({\bf R}^N)$, 
there exists a unique solution $v$ of 
$$
-\Delta v+v=a(x)v+f\quad\mbox{in}\quad{\bf R}^N,
\qquad
v\in H^1({\bf R}^N). 
$$ 
\end{lemma}
%
%%%%%%%%%%%%%%%%%%%%%%%%%%%%%%%%%%%%%
%%%%%%%%%%%%%%%%%%%%%%%%%%%%%%%%%%%%%
\section{Approximate solutions} \label{section:3}
%%%%%%%%%%%%%%%%%%%%%%%%%%%%%%%%%%%%%
%%%%%%%%%%%%%%%%%%%%%%%%%%%%%%%%%%%%%
Let $\mu$ be a nontrivial Radon measure in ${\bf R}^N$ with $\mbox{supp}\,\mu\subset B(0,R)$ for some $R>0$.  
Assume \eqref{eq:1.5}. 
Let $u$ be a $(C_{0} + L^{q}_{\rm c})$-solution of \eqref{eq:1.1}. 
Then 
\begin{equation}
\label{eq:3.1}
\mu_0(x):=[G*\mu](x)>0\quad\mbox{for almost all $x\in{\bf R}^N$}, \quad \mu_{0} \in C^\infty({\bf R}^{N} \setminus B(0, R)). 
\end{equation}
For any $\kappa>0$, 
we define $\{U_j^\kappa\}_{j=0}^\infty$ and $\{V_j^\kappa\}_{j=0}^\infty$ by 
\begin{equation}
\label{eq:3.2}
\begin{array}{ll}
U_0^\kappa:=\kappa\mu_0,\qquad
 & U^\kappa_j:=G*(U^\kappa_{j-1})^p+\kappa\mu_0,\qquad\,\,\,\, j=1,2,\dots,\vspace{3pt}\\
V_0^\kappa:=U_0^\kappa,\qquad
 & V^\kappa_j:=U^\kappa_j-U^\kappa_{j-1},\qquad\qquad\qquad j=1,2,\dots.
\end{array}
\end{equation}
By induction we easily obtain 
\begin{equation}
\label{eq:3.3}
\begin{split}
 & 0<U^\kappa_j(x)\le U^\kappa_{j+1}(x),\\
 & 0<U^\kappa_j(x)\le U^{\kappa'}_j(x)\quad\mbox{if}\quad \kappa\le\kappa',\\
 & U^{\kappa}_{j}(x) \le u(x), 
\end{split}
\end{equation}
for $j\in\{0,1,2,\dots\}$ and almost all $x\in{\bf R}^N$. 
In what follows, the letter $C$ denotes generic positive constants  
and it may have different values also within the same line. 
\begin{lemma}
\label{Lemma:3.1}
Assume the same conditions as in Theorem {\rm \ref{Theorem:1.1}}. 
\begin{itemize}
  \item[{\rm (i)}]
  For any $\kappa>0$, 
  there exists $c>0$ such that 
  $U_j^\kappa(x)\ge cg(x)$ for $j\in\{1,2,\dots\}$ and almost all $x\in {\bf R}^N$. 
  \item[{\rm (ii)}]
  For any $0<\kappa < \kappa'$ and $j\in\{0,1,2,\dots\}$, 
  \begin{equation}
  \label{eq:3.4}
  0<V^\kappa_j(x) < V^{\kappa'}_j(x)
  \end{equation} 
  for almost all $x\in {\bf R}^N$. 
\end{itemize}
\end{lemma}
{\bf Proof.}
We prove assertion~(i). 
Let $R>0$ be as in Theorem~\ref{Theorem:1.1}. 
Let $z\in {\bf R}^N\setminus B(0,R+2)$. 
It follows that 
$$
U_0^\kappa(x)=\kappa\int_{{\bf R}^N}G(x-y)\,d\mu(y)\ge \kappa\inf_{y\in B(0,R)}G(x-y)\,\int_{B(0,R)}\,d\mu\ge C
$$
for almost all $x\in B(z,1)$. 
Since $g=G*\chi_{B(0,1)}$, by \eqref{eq:2.2}, \eqref{eq:3.2} and \eqref{eq:3.3} we have
\begin{equation*}
\begin{split}
U^\kappa_j(x) & \ge U_1^\kappa(x)\ge \int_{{\bf R}^N} G(x-y)(U^\kappa_0(y))^p\,dy\\
 & \ge C\int_{B(z,1)} G(x-y)\,dy=Cg(x-z)\ge Cg(x)
\end{split}
\end{equation*}
for $j \in \{ 1, 2, \ldots\}$ and almost all $x\in{\bf R}^N$. 
Thus assertion~(i) follows. 

We prove assertion~(ii). 
Let $0 < \kappa < \kappa'$. 
It follows from \eqref{eq:3.1} and \eqref{eq:3.2} that  
\begin{align*} 
V_0^{\kappa'}-V_0^\kappa & =U_0^{\kappa'}-U_0^\kappa=(\kappa'-\kappa) \mu_0 >0,\\
V^{\kappa'}_{1} - V^{\kappa}_{1} & = ( U^{\kappa'}_{1} - U^{\kappa'}_{0}) - ( U^{\kappa}_{1} - U^{\kappa}_{0} )\\
 & = G * [ (U^{\kappa'}_{0})^{p} - (U^{\kappa'}_{0})^{p} ]
= G * [ (V^{\kappa'}_{0})^{p} - (V^{\kappa'}_{0})^{p} ] > 0,
\end{align*}
for almost all $x\in{\bf R}^N$. 
Then \eqref{eq:3.4} holds for $j=0, 1$. 

Assume that \eqref{eq:3.4} holds for some $j=j_0\in\{1,2,\dots\}$.
It follows from \eqref{eq:3.1} and \eqref{eq:3.2} that 
\begin{equation}
\label{eq:3.5}
\begin{split}
V^{\kappa'}_{j_0+1}-V^\kappa_{j_0+1}
 & =G*[(U^{\kappa'}_{j_0})^p-(U^{\kappa'}_{j_0-1})^p]-G*[(U^\kappa_{j_0})^p-(U^\kappa_{j_0-1})^p]\\
 & =G*[(U^{\kappa'}_{j_0-1}+V_{j_0}^{\kappa'})^p-(U^{\kappa'}_{j_0-1})^p]-G*[(U^\kappa_{j_0-1}+V_{j_0}^\kappa)^p-(U^\kappa_{j_0-1})^p].
\end{split}
\end{equation}
On the other hand, the function
$$
[0,\infty)\ni s\to(t+s)^p-s^p
$$
is monotone increasing for any fixed $t\ge 0$. 
Since \eqref{eq:3.4} holds for $j=j_0$, 
we deduce from \eqref{eq:3.3} and \eqref{eq:3.5} that 
\begin{equation*}
\begin{split}
V^{\kappa'}_{j_0+1}-V^\kappa_{j_0+1}
 & \ge G*[(U^\kappa_{j_0-1}+V_{j_0}^{\kappa'})^p-(U^\kappa_{j_0-1})^p]-G*[(U^\kappa_{j_0-1}+V_{j_0}^\kappa)^p-(U^\kappa_{j_0-1})^p]\\
 & =G*[(U^\kappa_{j_0-1}+V_{j_0}^{\kappa'})^p-(U^\kappa_{j_0-1}+V_{j_0}^\kappa)^p] > 0.
\end{split}
\end{equation*}
Thus \eqref{eq:3.4} holds for $j=j_0+1$. 
By induction we obtain \eqref{eq:3.4} for $j\in\{0,1,2,\dots\}$ and 
Lemma~\ref{Lemma:3.1} follows.
$\Box$\vspace{5pt}

Since 
$q>\max\{p,N(p-1)/2\}$,
we can find $r_*\in(1,\infty)$ such that 
\begin{equation}
\label{eq:3.6}
\max\left\{\frac{N}{2},\frac{q}{q-1}\right\}<r_*<\frac{q}{p-1},
\quad
\frac{1}{q}\not\in\left\{j\left(\frac{2}{N}-\frac{1}{r_*}\right)\,:\,j=0,1,2,\dots\right\}.
\end{equation}
Define a sequence $\{q_j\}_{j=0}^\infty$ by 
\begin{equation}
\label{eq:3.7}
\frac{1}{q_j}:=\frac{1}{q}-j\left(\frac{2}{N}-\frac{1}{r_*}\right).
\end{equation}
By \eqref{eq:3.6} there exists $j_*\in\{1,2,\dots\}$ such that 
\begin{equation}
\label{eq:3.8}
\frac{1}{q_{j_*-1}}>0>\frac{1}{q_{j_*}}.
\end{equation}
\begin{lemma}
\label{Lemma:3.2}
Assume the same condition as in Theorem {\rm \ref{Theorem:1.1}}. 
Let $K>0$. 
Then there exists $c_1>0$ such that 
\begin{equation}
\label{eq:3.9}
\begin{split}
 & \|V_j^\kappa\|_{L^{q_j}({\bf R}^N)}+\|g^{-1}V_j^\kappa\|_{L^\infty({\bf R}^N\setminus B(0,R))}\le c_1 \kappa^{(p-1)j+1},\\
 & \|U_j^\kappa\|_{L^q({\bf R}^N)}+\|g^{-1}U_j^\kappa\|_{L^\infty({\bf R}^N\setminus B(0,R))}\le c_1 \kappa,
\end{split}
\end{equation}
for $j=0,\dots,j_*-1$ and $0<\kappa<K$. 
Furthermore, for any $j\in\{j_*,j_*+1,\dots\}$, 
there exists $c_2>0$ such that 
\begin{equation} 
\label{eq:3.10}
\begin{split}
 & V_j^\kappa\in BC({\bf R}^N)\quad\mbox{and}\quad \|g^{-1}V_j^\kappa\|_{L^\infty({\bf R}^N)}\le c_2 \kappa^{(p-1)j+1},\\
 & \|U_j^\kappa\|_{L^q({\bf R}^N)}+\|g^{-1}U_j^\kappa\|_{L^\infty({\bf R}^N\setminus B(0,R))}\le c_2 \kappa,
\end{split}
\end{equation} 
for $0<\kappa<K$. 
\end{lemma}
{\bf Proof.}
Since $\mbox{supp}\,\mu$ is closed, 
we can find $R'\in(0,R)$ such that $\mbox{supp}\,\mu\subset B(0,R')$. 
Let $\{R_j\}_{j=0}^{j_*}$ be such that $R'<R_0<R_1<\cdots<R_{j_*}<R$. 
Then it follows from \eqref{eq:1.5} and \eqref{eq:2.2} that  
\begin{equation}
\label{eq:3.11}
\begin{split}
 & \|U_0^\kappa\|_{L^q({\bf R}^N)}=\kappa\|G*\mu\|_{L^q({\bf R}^N)}\le C\kappa,\\
 & 0<U_0^\kappa(x)\le\kappa\sup_{y\in B(0,R')}G(x-y)\int_{B(0,R')}\,d\mu\le C\kappa g(x)
 \quad\mbox{in}\quad {\bf R}^N\setminus B(0,R),
\end{split}
\end{equation}
for $0<\kappa<K$. 
Since $V_0^\kappa=U_0^\kappa$, 
by \eqref{eq:3.11} we have \eqref{eq:3.9} for $j=0$. 

Assume $j_*\ge 2$ and that 
\begin{equation}
\label{eq:3.12}
\begin{split}
 & \|V_j^\kappa\|_{L^{q_j}({\bf R}^N)}+\|g^{-1}V_j^\kappa\|_{L^\infty({\bf R}^N\setminus B(0,R_j))}\le C\kappa^{(p-1)j+1},\\
 & \|U_j^\kappa\|_{L^q({\bf R}^N)}+\|g^{-1}U_j^\kappa\|_{L^\infty({\bf R}^N\setminus B(0,R_j))}\le C\kappa,
\end{split}
\end{equation}
hold for some $j=j_0\in\{0,\dots,j_*-2\}$ and all $\kappa\in(0,K)$. 
It follows from \eqref{eq:3.2} and \eqref{eq:3.3} that
\begin{equation}
\label{eq:3.13}
\begin{split}
0<V_{j_0+1}^\kappa(x) & =U_{j_0+1}^\kappa(x)-U_{j_0}^\kappa(x)
= [G*\{(U_{j_0}^\kappa)^p-(U_{j_0-1}^\kappa)^p\}](x)\\
 & \le p [G*\{(U_{j_0}^\kappa)^{p-1}V_{j_0}^\kappa\}](x)\\
 & =p\biggr(\int_{B(0,R_{j_0})}+\int_{{\bf R}^N\setminus B(0,R_{j_0})}\biggr)\,G(x-y)U_{j_0}^\kappa(y)^{p-1}V_{j_0}^\kappa(y)\,dy
\end{split}
\end{equation}
for almost all $x\in{\bf R}^N$. 
Since $r_*<q/(p-1)$, 
we observe from \eqref{eq:3.12} that 
\begin{equation} \label{eq:3.14}
\begin{split}
\| (U_{j_0}^\kappa)^{p-1} \|_{L^{r_{*}}({\bf R}^{N})} 
& \le \| (U_{j_0}^\kappa)^{p-1} \|_{L^{r_{*}}(B(0,R_{j}))} + \| (U_{j_0}^\kappa)^{p-1} \|_{L^{r_{*}}({\bf R}^N\setminus B(0,R_j))} \\
& \le C \| U_{j_0}^\kappa \|_{L^{q}(B(0,R_{j}))}^{p-1} + C \kappa^{p-1} \| g^{p-1} \|_{L^{r_{*}}({\bf R}^N\setminus B(0,R_j))}
  \le C \kappa^{p-1}. 
\end{split}
\end{equation}
Then, by \eqref{eq:3.12} we have  
\begin{equation}
\label{eq:3.15}
\|(U_{j_0}^\kappa)^{p-1}V_{j_0}^\kappa\|_{L^{r_{j_0}}({\bf R}^N)}
\le\|(U_{j_0}^\kappa)^{p-1}\|_{L^{r_{*}}({\bf R}^N)}\|V_{j_0}^\kappa\|_{L^{q_{j_0}}({\bf R}^N)}
\le C\kappa^{(p-1)(j_0+1)+1}
\end{equation}
for $0<\kappa<K$, where $1/r_{j_0}=1/r_*+1/q_{j_0}$. 
Since 
\begin{equation} \label{eq:3.16}
\dfrac{1}{r_{j_{0}}} \le \dfrac{1}{r_{*}} + \dfrac{1}{q} < 1, \quad 
\frac{1}{r_{j_0}}-\frac{2}{N}
=\frac{1}{q_{j_0}}-\left(\frac{2}{N}-\frac{1}{r_*}\right)=\frac{1}{q_{j_0+1}}>0,
\end{equation}
by \eqref{eq:3.13} and \eqref{eq:3.15} we apply (G2) to obtain 
\begin{equation}
\label{eq:3.17}
\|V_{j_0+1}^\kappa\|_{L^{q_{j_0+1}}({\bf R}^N)}\le C\|(U_{j_0}^\kappa)^{p-1}V_{j_0}^\kappa\|_{L^{r_{j_0}}({\bf R}^N)}
\le C\kappa^{(p-1)(j_0+1)+1}
\end{equation}
for $0<\kappa<K$. 
This together with \eqref{eq:2.2} implies that 
\begin{equation}
\label{eq:3.18}
\begin{split}
 & \int_{B(0,R_{j_0})}\,G(x-y)U_{j_0}^\kappa(y)^{p-1}V_{j_0}^\kappa(y)\,dy\\
 & \le\sup_{y\in B(0,R_{j_0})}G(x-y)
\int_{B(0,R_{j_0})}\,U_{j_0}^\kappa(y)^{p-1}V_{j_0}^\kappa(y)\,dy\le C\kappa^{(p-1)(j_0+1)+1}g(x)
\end{split}
\end{equation}
for $x\in{\bf R}^N\setminus B(0,R_{j_{0}+1})$. 
Furthermore, by \eqref{eq:2.3} and \eqref{eq:3.12} 
we see that 
\begin{equation}
\label{eq:3.19}
\begin{split}
 & \int_{{\bf R}^N\setminus B(0,R_{j_0})}\,
 G(x-y)U_{j_0}^\kappa(y)^{p-1}V_{j_0}^\kappa(y)\,dy\\
 & \le C\kappa^{(p-1)(j_0+1)+1}[G*g^p](x)\le C\kappa^{(p-1)(j_0+1)+1}g(x),
 \quad x\in {\bf R}^N. 
\end{split}
\end{equation}
By \eqref{eq:3.13}, \eqref{eq:3.18} and \eqref{eq:3.19} we have 
\begin{equation}
\label{eq:3.20}
\|g^{-1}V_{j_0+1}^\kappa\|_{L^\infty({\bf R}^N\setminus B(0,R_{j_0+1}))}\le C\kappa^{(p-1)(j_0+1)+1}
\end{equation}
for $0<\kappa<K$. 
In addition, 
by \eqref{eq:3.2}, \eqref{eq:3.12}, \eqref{eq:3.17} and \eqref{eq:3.20} 
we obtain 
$$
\|U_{j_0+1}^\kappa\|_{L^q({\bf R}^N)}+\|g^{-1}U_{j_0+1}^\kappa\|_{L^\infty({\bf R}^N\setminus B(0,R_{j_0+1}))}\le C\kappa
$$
for $0<\kappa<K$. 
This together with \eqref{eq:3.17} and \eqref{eq:3.20} implies \eqref{eq:3.12} for $j=j_0+1$. 
Since \eqref{eq:3.9} holds for $j=0$, 
by induction we obtain \eqref{eq:3.12} for $j\in\{0,\dots,j_*-1\}$. 
This yields \eqref{eq:3.9} for $j\in\{0,\dots,j_*-1\}$. 

It remains to prove \eqref{eq:3.10}. 
Similarly to \eqref{eq:3.15}, 
by \eqref{eq:3.12} with $j=j_*-1$, 
we have
$$
\|(U_{j_*-1}^\kappa)^{p-1}V_{j_*-1}^\kappa\|_{L^{r_{j_*-1}}({\bf R}^N)}
\le\|(U_{j_*-1}^\kappa)^{p-1}\|_{L^{r_*}({\bf R}^N)}\|V_{j_*-1}^\kappa\|_{L^{q_{j_*-1}}({\bf R}^N)}
\le C\kappa^{(p-1)j_*+1}
$$
for $0<\kappa<K$, where $1/r_{j_*-1}=1/r_*+1/q_{j_*-1}<1$. 
On the other hand, 
it follows from \eqref{eq:3.7} and \eqref{eq:3.8} that  
$$
\frac{1}{r_{j*-1}}=\frac{1}{r_*}+\frac{1}{q}-(j_*-1)\left(\frac{2}{N}-\frac{1}{r_*}\right)
=\frac{2}{N}+\frac{1}{q_{j_*}}<\frac{2}{N}. 
$$ 
Then, similarly to \eqref{eq:3.13}, 
we see that 
\begin{equation}
\label{eq:3.21}
(U_{j_*-1}^\kappa)^p-(U_{j_*-2}^\kappa)^p\in L^{r_{j_*-1}}
\quad\mbox{with}\quad r_{j*-1}>\frac{N}{2}. 
\end{equation}
Combining \eqref{eq:3.2} with \eqref{eq:3.21}, we observe from (G3) that 
$$
V^\kappa_{j_*}\in C_0({\bf R}^N),
\qquad
\|V^\kappa_{j_*}\|_{L^\infty({\bf R}^N)}\le C\kappa^{(p-1)j_*+1}\quad\mbox{for}\quad 0<\kappa<K.
$$ 
Furthermore, 
similarly to \eqref{eq:3.20}, 
we obtain
$$
\|g^{-1}V_{j_*}^\kappa\|_{L^\infty({\bf R}^N\setminus B(0,R_{j_*}))}\le C\kappa^{(p-1)j_*+1}\quad\mbox{for}\quad 0<\kappa<K.
$$
This together with \eqref{eq:3.9} implies \eqref{eq:3.10} for $j=j_*$. 
Repeating this argument, we obtain \eqref{eq:3.10}. 
Thus Lemma~\ref{Lemma:3.2} follows.  
$\Box$

We obtain estimates on $U^{\kappa+\epsilon}_{j} - U^{\kappa}_{j}$ and $V^{\kappa+\epsilon}_{j} - V^{\kappa}_{j}$ in the following two lemmas. 
\begin{lemma}
\label{Lemma:3.3}
Assume the same condition as in Theorem {\rm \ref{Theorem:1.1}}. 
Then, for any $j\in\{0,1,2,\dots\}$ and $\kappa>0$, 
there exists $c > 0$ such that
\begin{equation}
\label{eq:3.22}
0\le U_j^{\kappa+\epsilon}-U_j^\kappa\le c \epsilon U_j^\kappa\quad\mbox{in}\quad{\bf R}^N
\quad\mbox{for}\quad 0<\epsilon\le 1.
\end{equation}
\end{lemma}
{\bf Proof.}
It follows from \eqref{eq:3.2} that 
$0\le U_0^{\kappa+\epsilon}-U_0^\kappa=\epsilon\mu_0=\kappa^{-1}\epsilon U_0^\kappa$ in ${\bf R}^N$. 
So \eqref{eq:3.22} holds for $j=0$ and $0<\epsilon\le 1$. 

Assume that \eqref{eq:3.22} holds for some $j=j_0\in\{0,1,2,\dots\}$ and all $0<\epsilon\le 1$. 
Similarly to \eqref{eq:3.13}, by \eqref{eq:3.2} we have 
\begin{equation*} 
\begin{split}
 U_{j_0+1}^{\kappa+\epsilon}-U_{j_0+1}^\kappa
 & =G*[(U_{j_0}^{\kappa+\epsilon})^p-(U_{j_0}^\kappa)^p]+\epsilon\mu_0\\
 & \le pG*[(U_{j_0}^{\kappa+\epsilon})^{p-1}(U_{j_0}^{\kappa+\epsilon}-U_{j_0}^\kappa)]+\epsilon\mu_0\\
 & \le C G * [((1+\epsilon)U^{\kappa}_{j_{0}})^{p-1} \epsilon U^{\kappa}_{j_{0}}] + \kappa^{-1} \epsilon U^{\kappa}_{0} 
 \le C\epsilon U_{j_0+1}^\kappa\quad\mbox{in}\quad{\bf R}^N
\end{split}
\end{equation*}
for $0<\epsilon\le 1$. 
This implies that \eqref{eq:3.22} holds for $j=j_0+1$ and $0<\epsilon\le 1$.
Therefore, by induction \eqref{eq:3.22} holds for $j\in\{0,1,2,\dots\}$ and $0<\epsilon\le 1$. 
Thus Lemma~\ref{Lemma:3.3} follows.
$\Box$
\begin{lemma}
\label{Lemma:3.4}
Assume the same conditions as in Theorem~{\rm\ref{Theorem:1.1}}.
Let $\kappa>0$. 
Then there exists $c_1>0$ such that 
\begin{equation}
\label{eq:3.23}
\left\|V_j^{\kappa+\epsilon}-V_j^\kappa\right\|_{L^{q_j}({\bf R}^N)}
+\left\|g^{-1}[V_j^{\kappa+\epsilon}-V_j^\kappa]\right\|_{L^\infty({\bf R}^N\setminus B(0,R))}
\le c_1\epsilon^{\min\{1,p-1\}},
\end{equation}
for $j=0,\dots,j_*-1$ and $0<\epsilon\le 1$. 
Furthermore, for any $j\in\{j_*,j_*+1,\dots\}$, 
there exists $c_2>0$ such that  
\begin{equation} 
\label{eq:3.24}
\left\|g^{-1}[V_j^{\kappa+\epsilon}-V_j^\kappa]\right\|_{L^\infty({\bf R}^N)}\le c_2\epsilon^{\min\{1,p-1\}}
\end{equation} 
for $0<\epsilon\le 1$. 
\end{lemma}
{\bf Proof.}
Let $\{R_j\}_{j=0}^{j_*}$ and $\{r_j\}$ be as in the proof of Lemma~{\rm\ref{Lemma:3.2}}. 
In the case of $j=0$, 
since $V_j^{\kappa+\epsilon}-V_j^\kappa=\kappa^{-1}\epsilon U^{\kappa}_0$, 
by Lemma~\ref{Lemma:3.2} we have \eqref{eq:3.23} for $0 \le \epsilon \le 1$. 

Let $j_*\ge 1$ and $j=1$. 
By \eqref{eq:3.2} we see that 
\begin{align*}
V^{\kappa+\epsilon}_{1} - V^{\kappa}_{1} 
 &= G * [ (U^{\kappa + \epsilon}_{0})^{p} - (U^{\kappa}_{0})^{p} ] \\
 &= G * [ p (U^{\kappa+\epsilon}_{0})^{p-1} (U^{\kappa+\epsilon}_{0} - U^{\kappa}_{0} ) ]  \\
 &\le C \epsilon G * [ (U^{\kappa}_{0})^{p} ] 
  = C \epsilon ( U^{\kappa}_{1} - U^{\kappa}_{0} ) 
   = C \epsilon V^{\kappa}_{0}. 
\end{align*}
Then Lemma~\ref{Lemma:3.2} implies \eqref{eq:3.23} with $j=1$. 

Let $j_*\ge 2$ and assume that 
\begin{equation}
\label{eq:3.25}
\left\|V_j^{\kappa+\epsilon}-V_j^\kappa\right\|_{L^{q_j}({\bf R}^N)}
+\left\|g^{-1}[V_j^{\kappa+\epsilon}-V_j^\kappa]\right\|_{L^\infty({\bf R}^N\setminus B(0,R_j))}
\le C\epsilon^{\min\{1,p-1\}}
\end{equation}
for some $j=j_0\in\{0,\dots,j_*-2\}$ and all $0<\epsilon\le 1$. 
It follows from \eqref{eq:3.2} that
\begin{equation}
\label{eq:3.26}
\begin{split}
V_{j_0+1}^{\kappa+\epsilon}-V_{j_0+1}^\kappa
 & =G*[(U_{j_0}^{\kappa+\epsilon})^p-(U_{j_0-1}^{\kappa+\epsilon})^p]
 -G*[(U_{j_0}^\kappa)^p-(U_{j_0-1}^\kappa)^p]\\
 & =G*[h(1)-h(0)],
\end{split}
\end{equation}
where
\begin{equation*}
\begin{split}
 & h(t):=\alpha(t)^p-\beta(t)^p,\\
 & \alpha(t):=tU_{j_0}^{\kappa+\epsilon}+(1-t)U_{j_0-1}^{\kappa+\epsilon},\qquad
\beta(t):=tU_{j_0}^\kappa+(1-t)U_{j_0-1}^\kappa.
\end{split}
\end{equation*}
On the other hand,
\begin{equation}
\label{eq:3.27}
\begin{split}
h'(t) & =p\alpha(t)^{p-1}V_{j_0}^{\kappa+\epsilon}-p\beta(y)^{p-1}V_{j_0}^\kappa\\
 & =p\alpha(t)^{p-1}(V_{j_0}^{\kappa+\epsilon}-V_{j_0}^\kappa)
+p(\alpha(t)^{p-1}-\beta(t)^{p-1})V_{j_0}^\kappa.
\end{split}
\end{equation}
It follows from Lemma~\ref{Lemma:3.3}, \eqref{eq:3.3} and \eqref{eq:3.14} that 
\begin{equation} \label{eq:3.28}
0\le \alpha(t)^{p-1}\le (U_{j_0}^{\kappa+\epsilon}+U_{j_0-1}^{\kappa+\epsilon})^{p-1}
\le C(U_{j_0}^{\kappa+\epsilon})^{p-1}\le C(U_{j_{0}}^\kappa)^{p-1}\in L^{r_*}({\bf R}^N).
\end{equation}
Then, similarly to \eqref{eq:3.15}, 
by \eqref{eq:3.25} we obtain
\begin{equation}
\label{eq:3.29}
\begin{split}
 & \|\alpha(t)^{p-1}(V_{j_0}^{\kappa+\epsilon}-V_{j_0}^\kappa)\|_{L^{r_{j_0}}({\bf R}^N)}\\
 & \le\|\alpha(t)^{p-1}\|_{L^{r_*}({\bf R}^N)}\|V_{j_0}^{\kappa+\epsilon}-V_{j_0}^\kappa\|_{L^{q_{j_0}}({\bf R}^N)}
\le C\epsilon^{\min\{1,p-1\}}
\end{split}
\end{equation}
for $0<\epsilon\le 1$. 
On the other hand, by Lemma~\ref{Lemma:3.3} and \eqref{eq:3.3} we see that  
\begin{equation*}
\begin{split}
0\le\alpha(t)-\beta(t) & =t(U_{j_0}^{\kappa+\epsilon}-U_{j_0}^\kappa)
+(1-t)(U_{j_0-1}^{\kappa+\epsilon}-U_{j_0-1}^\kappa)\\
 & \le C\epsilon(tU_{j_0}^\kappa+(1-t)U_{j_0-1}^\kappa)
\le C\epsilon U_{j_0}^\kappa,\\
\max\{\alpha(t),\beta(t)\} & =\alpha(t)\le U_{j_0}^{\kappa+\epsilon}\le CU_{j_0}^\kappa,
\end{split}
\end{equation*}
for $0<\epsilon\le 1$. 
In the case of $1<p\le 2$, we have
\begin{equation}
\label{eq:3.30}
\begin{split}
\alpha(t)^{p-1}-\beta(t)^{p-1}
 & =(\alpha(t)-\beta(t)+\beta(t))^{p-1}-\beta(t)^{p-1}\\
 & \le(\alpha(t)-\beta(t))^{p-1}
 \le C\epsilon^{p-1}(U_{j_0}^\kappa)^{p-1}.
\end{split}
\end{equation}
On the other hand, in the case of $p>2$, 
\begin{equation}
\label{eq:3.31}
\alpha(t)^{p-1}-\beta(t)^{p-1}\le(p-1)\alpha(t)^{p-2}(\alpha(t)-\beta(t))
\le C\epsilon(U_{j_0}^\kappa)^{p-1}.
\end{equation}
Similarly to \eqref{eq:3.14}, 
by \eqref{eq:3.9}, \eqref{eq:3.30} and \eqref{eq:3.31} we see that 
\begin{equation}
\label{eq:3.32}
\begin{split}
 & \|(\alpha(t)^{p-1}-\beta(t)^{p-1})V_{j_0}^\kappa\|_{L^{r_{j_0}}({\bf R}^N)}\\
 & \le\|(\alpha(t)^{p-1}-\beta(t)^{p-1})\|_{L^{r_*}({\bf R}^N)}\|V_j^\kappa\|_{L^{q_{j_0}}({\bf R}^N)}
 \le C\epsilon^{\min\{p-1,1\}}
\end{split}
\end{equation}
for $0<\epsilon\le 1$. 
Therefore, by \eqref{eq:3.27}, \eqref{eq:3.29} and \eqref{eq:3.32} 
we have 
\begin{equation}
\label{eq:3.33}
\|h(1)-h(0)\|_{L^{r_{j_0}}({\bf R}^N)}\le\int_0^1\|h'(t)\|_{L^{r_{j_{0}}}({\bf R}^N)}\,dt\le C\epsilon^{\min\{1,p-1\}}
\end{equation}
for $0<\epsilon\le 1$. 
Then, similarly to \eqref{eq:3.17}, 
we deduce from (G2), \eqref{eq:3.16} and \eqref{eq:3.26} that 
\begin{equation} \label{eq:3.34}
\left\|V_{j_0+1}^{\kappa+\epsilon}-V_{j_0+1}^\kappa\right\|_{L^{q_{j_0+1}}({\bf R}^N)}
\le C\|h(1)-h(0)\|_{L^{r_{j_0}}({\bf R}^N)}
\le C\epsilon^{\min\{1,p-1\}}
\end{equation}
for $0<\epsilon\le 1$.

On the other hand, it follows that 
$$
0\le V_{j_0+1}^{\kappa+\epsilon}-V_{j_0+1}^\kappa=\left(\int_{B(0,R_{j_0})}+\int_{{\bf R}^N\setminus B(0,R_{j_0})}\right)
G(x-y)(h(y,1)-h(y,0))\,dy.
$$
By \eqref{eq:3.33} we have 
\begin{equation}
\label{eq:3.35}
\begin{split}
 & \int_{B(0,R_{j_{0}})}G(x-y)(h(y,1)-h(y,0))\,dy\\
 & \le\sup_{y\in B(0,R_{j_0})}G(x-y)\int_{B(0,R_{j_0})}(h(y,1)-h(y,0))\,dy\le Cg(x)
\end{split}
\end{equation}
for $x\in{\bf R}^N\setminus B(0,R_{j_0+1})$. 
Since \eqref{eq:3.12} holds for $j=j_0$, 
from \eqref{eq:3.25}, \eqref{eq:3.28}, \eqref{eq:3.30} and \eqref{eq:3.31} we deduce that
\begin{equation*}
\begin{split}
 & 0\le \alpha(t)^{p-1}(V_{j_0}^{\kappa+\epsilon}-V_{j_0}^\kappa)\le C\epsilon^{\min\{1,p-1\}}g^p,\\
 & 0\le (\alpha(t)^{p-1}-\beta(t)^{p-1})V_{j_0}^\kappa\le C\epsilon^{\min\{1,p-1\}}g^p,
\end{split}
\end{equation*}
for $x\in{\bf R}^N\setminus B(0,R_{j_0})$.
This together with \eqref{eq:3.27} implies that 
\begin{equation}
\label{eq:3.36}
\int_{{\bf R}^N\setminus B(0,R_{j_0})}G(x-y)(h(y,1)-h(y,0))\,dy
\le C \epsilon^{\min\{1,p-1\}} G*g^p\le C\epsilon^{\min\{1,p-1\}} g.
\end{equation}
By \eqref{eq:3.35} and \eqref{eq:3.36} we obtain 
\begin{equation} \label{eq:3.37}
\left\|g^{-1}[V_{j_0+1}^{\kappa+\epsilon}-V_{j_0+1}^\kappa]\right\|_{L^\infty({\bf R}^N\setminus B(0,R_{j_0 +1}))}
\le C\epsilon^{\min\{1,p-1\}}
\end{equation}
for $0<\epsilon\le 1$. 
This together with \eqref{eq:3.34} implies that \eqref{eq:3.25} holds for $j=j_0+1$. 
Therefore, by induction we see that \eqref{eq:3.25} holds for $j\in\{0,1,\dots,j_*-1\}$ and $0<\epsilon\le 1$. 

We prove \eqref{eq:3.24}. 
By \eqref{eq:3.25} with $j=j_*-1$, 
similarly to \eqref{eq:3.33}, we have 
$$
\|h(1)-h(0)\|_{L^{r_{j_*-1}}({\bf R}^N)}\le\int_0^1\|h'(t)\|_{L^{r_j}({\bf R}^N)}\,dt\le C\epsilon^{\min\{1,p-1\}}
$$
for $0<\epsilon\le 1$. 
It follows from \eqref{eq:3.21} that $r_{j_*-1}>N/2$. 
Then we deduce from (G2) and \eqref{eq:3.26} that 
\begin{equation}
\label{eq:3.38}
\left\|V_{j_0+1}^{\kappa+\epsilon}-V_{j_0+1}^\kappa\right\|_{L^\infty({\bf R}^N)}\le C\epsilon^{\min\{1,p-1\}}
\end{equation}
for $0<\epsilon\le 1$.
Furthermore, similarly to \eqref{eq:3.37}, we have 
\begin{equation}
\label{eq:3.39}
\left\|g^{-1}[V_{j_*}^{\kappa+\epsilon}-V_{j_*}^\kappa]\right\|_{L^\infty({\bf R}^N\setminus B(0,R_{j_*}))}
\le C\epsilon^{\min\{1,p-1\}}
\end{equation}
for $0<\epsilon\le 1$. 
By \eqref{eq:3.38} and \eqref{eq:3.39} we have \eqref{eq:3.24}. 
Thus Lemma~\ref{Lemma:3.4} follows.
$\Box$
%%%%%%%%%%%%%%%%%%%%%%%%%%%%%%%%%%%%%
%%%%%%%%%%%%%%%%%%%%%%%%%%%%%%%%%%%%%
\section{Proof of Theorem~\ref{Theorem:1.1}} \label{section:4}
%%%%%%%%%%%%%%%%%%%%%%%%%%%%%%%%%%%%%
%%%%%%%%%%%%%%%%%%%%%%%%%%%%%%%%%%%%%
For any $\kappa>0$, define
$$
S_\kappa:=\{u\,:\,\mbox{$u$ is a $(C_0+L^q_{\rm c})$-solution of \eqref{eq:1.1}}\}.
$$
According to Definition \ref{Definition:1.1}, \eqref{eq:3.2} and \eqref{eq:3.3}, 
the following lemma holds. 
\begin{lemma}
\label{Lemma:4.1}
Assume the same condition as in Theorem {\rm \ref{Theorem:1.1}}. Let $\kappa>0$. 
Then the following conditions are equivalent{\rm :}
\begin{itemize}
  \item[{\rm (i)}] 
  $u=w+U^\kappa_{j_*}$ is a $(C_0+L^q_{\rm c})$-solution of \eqref{eq:1.1}{\rm ;}
  \item[{\rm (ii)}]
  $w\in C_0({\bf R}^N)$ is positive in ${\bf R}^N$ and $w$ satisfies 
  \begin{equation}
  \label{eq:4.1}
  w=G*[(w+U_{j_*}^\kappa)^p-(U_{j_*-1}^\kappa)^p]\quad\mbox{in}\quad{\bf R}^N. 
  \end{equation}
\end{itemize}
\end{lemma}
Furthermore, we have:
\begin{lemma}
\label{Lemma:4.2}
Assume the same condition as in Theorem {\rm \ref{Theorem:1.1}}. Let $\kappa>0$. 
If $w\in C_0({\bf R}^N)$ is a positive solution of \eqref{eq:4.1}, then 
there exists $c>0$ such that 
\begin{align} \label{eq:4.2}
0\le w(x)\le c g(x)\quad\mbox{in}\quad{\bf R}^N.
\end{align}
Furthermore, 
$w\in H^1({\bf R}^N)$ and $w$ is a weak solution of  
\begin{equation}
\label{eq:4.3}
-\Delta w+w=(w+U_{j_*}^\kappa)_+^p-(U_{j_*-1}^\kappa)^p\quad\mbox{in}\quad{\bf R}^N,
\end{equation}
that is
\begin{equation}
\label{eq:4.4}
\int_{{\bf R}^N}[\nabla w\cdot\nabla\psi+w\psi]\,dx
=\int_{{\bf R}^N}[(w+U_{j_*}^\kappa)_+^p-(U_{j_*-1}^\kappa)^p]\psi\,dx
\end{equation}
for $\psi\in H^1({\bf R}^N)$. Here $s_{+}:= \max\{ s, 0 \}$ for $s \in {\bf R}$. 
\end{lemma}
{\bf Proof.}
Let $R>0$ be as in Theorem \ref{Theorem:1.1}. 
Let $\zeta\in C^\infty({\bf R}^N)$ be such that 
$$
0\le \zeta\le 1\quad\mbox{in}\quad{\bf R}^N,\qquad
\zeta=1\quad\mbox{in}\quad B(0,R+1),\qquad
\zeta=0\quad\mbox{in}\quad{\bf R}^N\setminus B(0,R+2).
$$
Let $W:=(w+U_{j_*}^\kappa)^p-(U_{j_*-1}^\kappa)^p$. Set 
\begin{equation*}
w_1(x):=\int_{{\bf R}^N}G(x-y) \zeta(y) W(y) \,dy,
\qquad
w_2(x):=\int_{{\bf R}^N}G(x-y)(1-\zeta(y))W(y)\,dy.
\end{equation*}
Since 
$$
W \le p (w + U^{\kappa}_{j_{*}})^{p-1} (w + V^{\kappa}_{j_{*}}),
$$ 
it follows from \eqref{eq:3.9} and \eqref{eq:3.10} that $\zeta W \in L^1({\bf R}^N)\cap L^r({\bf R}^N)$ for some $r>N/2$. 
Thus we deduce from (G4) that 
\begin{equation}
\label{eq:4.5}
w_1\in C_0({\bf R}^N)\cap L^1({\bf R}^N)\cap H^1({\bf R}^N).
\end{equation}
Furthermore, 
\begin{equation}
\label{eq:4.6}
w_1(x)\le \sup_{y\in B(0,R+2)}G(x-y)\int_{B(0,R+2)}W(y)\,dy
\le Cg(x) 
\end{equation}
for $x\in {\bf R}^N\setminus B(0,R+3)$. 

On the other hand, by \eqref{eq:2.1}, \eqref{eq:3.1}, \eqref{eq:3.2} and Lemma \ref{Lemma:3.2}, 
we deduce from a bootstrap argument that $w_2\in C^{2,\theta}({\bf R}^N)$, where $0 < \theta < 1$, and $w_2$ is a classical solution of 
%we see that $w_2 \in C^1({\bf R}^N)$. 
%Since it follows from \eqref{eq:3.1}, \eqref{eq:3.2} and Lemma \ref{Lemma:3.2} that $(1-\zeta) W \in C^1({\bf R}^N)$, 
%we deduce from \eqref{eq:2.1} that $w_2\in C^{2,\theta}({\bf R}^N)$, where $0 < \theta < 1$, and $w_2$ is a classical solution of  
$$
-\Delta w_2+w_2=(1-\zeta)W \quad \mbox{in} \quad {\bf R}^N.
$$
Let $x_0\in {\bf R}^N\setminus B(0,R+3)$. 
Since $w \in C_{0}({\bf R}^{N})$, by Lemma \ref{Lemma:3.1} (i), \eqref{eq:3.1}, \eqref{eq:3.3} and \eqref{eq:3.10} 
we can find a positive constant $C>0$ such that 
\begin{align*}
W \ge p (U^{\kappa}_{j_{*}-1})^{p-1}(w + V^{\kappa}_{j_{*}}) \ge C \quad \text{in} \quad B(x_0,1). 
\end{align*}
Thus we have 
\begin{equation*}
\begin{split}
w_2(x) & \ge\int_{B(x_0,1)}G(x-y)W(y)\,dy\\
 & \ge C\int_{B(x_0,1)}G(x-y)\,dy=Cg(x-x_0)\ge Cg(x),\qquad x\in{\bf R}^N.
\end{split}
\end{equation*}
Since $w=w_1+w_2$, it follows from \eqref{eq:3.10} and \eqref{eq:4.6} that
\begin{equation} \label{eq:4.7}
\begin{split}
\frac{(1-\zeta)W}{w_2}
\le p (w+U_{j_*}^\kappa)^{p-1}\frac{w+V_{j_*}^\kappa}{w_2} 
\le C p (w+U_{j_*}^\kappa)^{p-1} \quad \text{in} \quad {\bf R}^N\setminus B(0,R+3). 
\end{split}
\end{equation}
Let $\epsilon\in(0,1)$. Since $w \in C_{0}({\bf R})$, by \eqref{eq:3.10} and \eqref{eq:4.7} we can find $L>R+3$ such that 
$$
-\Delta w_2+w_2=(1-\zeta)W
\le C p (w+U_{j_*}^\kappa)^{p-1}w_2\le\epsilon w_2\quad \text{in} \quad {\bf R}^N\setminus B(0,L). 
$$
Let $\gamma>0$ be such that 
$$
\gamma G(\sqrt{1-\epsilon}x)\ge w_2(x) \quad\mbox{on}\quad\partial B(0,L). 
$$
Set 
$$
z(x):=w_2(x)-\gamma G(\sqrt{1-\epsilon}x) \quad \text{for} \quad x\in{\bf R}^N\setminus B(0,L).
$$
Since  
$$
-\Delta z + (1-\epsilon)z \le 0\quad\mbox{in}\quad{\bf R}^N\setminus\overline{B(0,L)},
\quad
z\le 0\quad\mbox{on}\quad\partial B(0,L),
\quad
\lim_{x\to\infty}z(x)=0,
$$
we deduce from the maximum principle that $z\le 0$ in ${\bf R}^N\setminus B(0,L)$. 
Then we have 
$$
0<w_2(x)\le \gamma G(\sqrt{1-\epsilon}x)\le Ce^{-{\sqrt{1-\epsilon}|x|}}\quad\mbox{in}\quad{\bf R}^N \setminus B(0,L).
$$
Since $\epsilon$ is arbitrary and $p>1$, 
we can find $\sigma>1$ such that 
$$
0\le W = [(w_1+w_2+U_{j_*}^\kappa)^p-(U_{j_*-1}^\kappa)^p]\le Cg^\sigma
\quad\mbox{in}\quad{\bf R}^N \setminus B(0,L).
$$
Then it follows from \eqref{eq:2.3} that  
\begin{align*}
w(x) &= \int_{B(0,L)} G(x-y) W(y) \, dy + \int_{{\bf R}^{N} \setminus B(0,L)} G(x-y) W(y) \, dy \\
 &\le \sup_{y \in B(0,L)} G(x-y) \int_{B(0,L)} W(y) \, dy 
        + \int_{{\bf R}^{N} \setminus B(0,L)} G(x-y) g(y)^{\sigma} \, dy 
 \le C g(x) 
\end{align*}
for $x \in {\bf R}^{N} \setminus B(0,L+1)$. 
Since $w$ is continuous in ${\bf R}^{N}$, we obtain \eqref{eq:4.2}. 
Moreover, we observe that $W \in L^{1}({\bf R}^{N}) \cap L^{r}({\bf R}^{N})$ for some $r > N/2$. 
Then (G4) implies that $w\in H^1({\bf R}^N)$. 
Furthermore, 
\begin{equation*}
\begin{split}
\int_{{\bf R}^N}[\nabla w\nabla\psi+w\psi]\,dx
& =\int_{{\bf R}^N}[-w\Delta \psi+w\psi]\,dx\\
 & =\int_{{\bf R}^N}\int_{{\bf R}^N}
G(x-y)W(y)[-\Delta_x\psi(x)+\psi(x)]\,dx\,dy\\
 & =\int_{{\bf R}^N}W(y)\psi(y)\,dy
 \quad \text{for} \quad \psi \in C_{\rm c}^\infty({\bf R}^N).
\end{split}
\end{equation*}
Therefore $w$ is a weak solution of \eqref{eq:4.3} in ${\bf R}^N$. 
Thus Lemma~\ref{Lemma:4.2} follows.
$\Box$
\vspace{5pt}

We show that the existence of supersolutions of \eqref{eq:1.1} ensures
the existence of the solution of \eqref{eq:1.1}. 
\begin{lemma}
\label{Lemma:4.3}
Assume the same condition as in Theorem {\rm \ref{Theorem:1.1}}. 
Let $v$ be a $(C_0+L^q_{\rm c})$-supersolution of \eqref{eq:1.1}. 
Then 
\begin{equation}
\label{eq:4.8}
U_j^\kappa(x)\le v(x)
\end{equation}
for almost all $x\in{\bf R}^N$ and all $j\in\{0,1,2,\dots\}$. 
Furthermore, 
$$
U^\kappa_\infty(x):=\lim_{j\to\infty}U_j^\kappa(x)
$$
exists for almost all $x\in{\bf R}^N$ and it is  a minimal $(C_0+L^q_{\rm c})$-solution of \eqref{eq:1.1}.
\end{lemma}
{\bf Proof.}
Let $R$ be as in Theorem \ref{Theorem:1.1}. 
Since $v$ is a $(C_0+L^q_{\rm c})$-supersolution of \eqref{eq:1.1}, we can find $L>R$ such that 
\begin{equation}
\label{eq:4.9}
v\in C({\bf R}^N\setminus B(0,L)),
\qquad
v\in L^q(B(0,L)).
\end{equation}
Similarly to \eqref{eq:3.3}, by induction we obtain \eqref{eq:4.8}. 
Furthermore, by \eqref{eq:3.2}, \eqref{eq:3.3} and \eqref{eq:4.8} 
we see that $U^\kappa_\infty$ exists and satisfies
\begin{equation}
\label{eq:4.10}
U^\kappa_\infty(x)=[G*(U^\kappa_\infty)^p](x)+\kappa [G*\mu](x),
\qquad
0<U^\kappa_\infty(x)\le v(x), 
\end{equation}
for almost all $x\in {\bf R}^N$. 
Setting $w_\infty = U^\kappa_\infty - U^\kappa_{j_*}$, we have 
\begin{align*}
w_\infty = G*[(w_\infty + U^\kappa_{j_*})^p - (U^\kappa_{j_{*}-1})^p] \quad \text{in} \quad {\bf R}^N. 
\end{align*}
Let $\zeta\in C^\infty({\bf R}^N)$ be such that 
$$
0\le \zeta\le 1\quad\mbox{in}\quad{\bf R}^N,\qquad
\zeta=1\quad\mbox{in}\quad B(0,L+1),\qquad
\zeta=0\quad\mbox{in}\quad{\bf R}^N\setminus B(0,L+2).
$$
Let $W_\infty = (w_\infty + U^\kappa_{j_*})^p - (U^\kappa_{j_{*}-1})^p$ and set 
$$
I(x):=\int_{{\bf R}^N} G(x-y) \zeta(y) W_\infty(y) \,dy,
\quad
J(x):=\int_{{\bf R}^N}G(x-y) (1-\zeta(y)) W_\infty(y) \,dy.
$$
By \eqref{eq:4.9} and \eqref{eq:4.10} we see that $U^\kappa_\infty\in L^q(B(0,L+2))$. 
Similarly to \eqref{eq:4.5}, it follows from (G3) and $q > N/2$ that $I\in C_0({\bf R}^N)$. 
Since $v\in C_{0}({\bf R}^{N}) + L^{q}_{{\rm c}}({\bf R}^{N})$, 
by Lemma~\ref{Lemma:3.2} and \eqref{eq:4.10} we have $(1-\zeta)W_\infty \in L^\infty({\bf R}^N)$. 
Then we deduce from \eqref{eq:1.4} that $J \in C({\bf R}^N)$. 
Therefore we obtain $U^\kappa_\infty \in C({\bf R}^N \setminus B(0, L+2))$. 
Combining $v\in C_{0}({\bf R}^{N}) + L^{q}_{{\rm c}}({\bf R}^{N})$ with \eqref{eq:4.10}, we observe that $U^\kappa_\infty(x) \to 0$ as $x \to \infty$. 
Thus $U_\infty^\kappa$ is a minimal $(C_0+L^q_{\rm c})$-solution of \eqref{eq:1.1}, 
and the proof is complete.
$\Box$
\vspace{5pt}

By the use of Lemma~\ref{Lemma:4.3} we prove the following lemma. 
\begin{lemma}
\label{Lemma:4.4}
Assume the same conditions as in Theorem~{\rm\ref{Theorem:1.1}}. 
Then
\begin{itemize}
  \item[{\rm (i)}] 
  $S_\kappa\not=\emptyset$ for sufficiently small $\kappa>0${\rm ;}
  \item[{\rm (ii)}] 
  If $S_\kappa\not=\emptyset$, then problem~\eqref{eq:1.1} possesses a minimal $(C_0+L^q_{\rm c})$-solution $u^{\kappa}$.
\end{itemize}
\end{lemma}
{\bf Proof.}
We prove assertion~(i). Let 
$$
gBC({\bf R}^N):=\{f\,:\,g^{-1}f\in BC({\bf R}^N)\}. 
$$
Then $gBC({\bf R}^{N})$ is a Banach space with the norm $|||f|||:=\sup_{x\in{\bf R}^N}|g(x)^{-1}f(x)|$. 
Let $0<\kappa<1$ and define
$$
F_\kappa[v]:=G*[(v+U_{j_*}^\kappa)_+^p-(U_{j_*-1}^\kappa)^p]
\quad \text{for} \quad v\in gBC({\bf R}^N).
$$
For any $v\in gBC({\bf R}^N)$, we have
\begin{equation}
\label{eq:4.11}
0\le (|v|+U_{j_*}^\kappa)^p-(U_{j_*-1}^\kappa)^p
\le p (|v|+U_{j_*}^\kappa)^{p-1}(|v|+V_{j_*}^\kappa). 
\end{equation}
It follows from Lemma~\ref{Lemma:3.2} that
\begin{equation}
\label{eq:4.12}
\||v|+U_{j_*}^\kappa\|_{L^q({\bf R}^N)}\le C|||v|||+C\kappa,
\quad
||||v|+V_{j_*}^\kappa|||\le |||v|||+C\kappa^{(p-1)j_*+1},
\end{equation}
for $0<\kappa<1$. 
Thus, by \eqref{eq:4.11} and \eqref{eq:4.12} we have $(v+U_{j_*}^\kappa)_+^p-(U_{j_*-1}^\kappa)^p \in L^{q/(p-1)}({\bf R}^N)$. 
Then it follows from (G3) and $q/(p-1) > N/2$ that $F_\kappa[v]\in C_{0}({\bf R}^N)$. 
Set  
\begin{align*}
I_{1}(x) &:= \int_{B(0,R+1)} G(x-y) [(v+U_{j_*}^\kappa)_+^p-(U_{j_*-1}^\kappa)^p](y) \, dy, \\ 
I_{2}(x) &:= \int_{{\bf R}^{N} \setminus B(0,R+1)} G(x-y) [(v+U_{j_*}^\kappa)_+^p-(U_{j_*-1}^\kappa)^p](y) \, dy. 
\end{align*}
By \eqref{eq:4.11} and \eqref{eq:4.12} we have 
\begin{align*}
I_{1}(x) & \le p \int_{B(0, R+1)} G(x-y) (|v|+U_{j_*}^\kappa)^{p-1}(y) g(y) (g^{-1} |v| + g^{-1} V_{j_*}^\kappa)(y) \, dy \\
 & \le C \||v|+U_{j_*}^\kappa\|_{L^q({\bf R}^N)}^{p-1} ||| |v| + V_{j_*}^\kappa||| \\
 & \le C \||v|+U_{j_*}^\kappa\|_{L^q({\bf R}^N)}^{p-1} ||| |v| + V_{j_*}^\kappa||| g(x)
\end{align*}
for $x \in B(0, R+2)$. 
Similarly, it follows that  
\begin{align*}
I_{1}(x) & \le p \sup_{y \in B(0, R+1)} G(x-y) \int_{B(0, R+1)} (|v|+U_{j_*}^\kappa)^{p-1} g g^{-1} (|v|+V_{j_*}^\kappa) \, dy \\
 & \le C \| |v|+U_{j_*}^\kappa \|^{p-1}_{L^{q}({\bf R}^{N})} ||| |v|+V_{j_*}^\kappa ||| g(x) \\
 & \le C (|||v|||+C\kappa)^{p-1}(|||v|||+C\kappa^{(p-1)j_*+1}) g(x)
\end{align*}
for $x \in {\bf R}^{N} \setminus B(0, R+2)$. 
Similarly, we observe from \eqref{eq:2.3} and Lemma \ref{Lemma:3.2} that  
\begin{align*}
I_{2}(x) & \le p \int_{{\bf R}^{N} \setminus B(0,R+1)} G(x-y) g(y)^{p} (g^{-1}|v|+ g^{-1}U_{j_*}^\kappa)^{p-1}(y) (g^{-1}|v|+ g^{-1}V_{j_*}^\kappa)(y) \, dy \\
 & \le C (|||v|||+C\kappa)^{p-1}(|||v|||+C\kappa^{(p-1)j_*+1}) \int_{{\bf R}^{N} \setminus B(0,R+1)} G(x-y) g(y)^{p} \, dy \\
 & \le C (|||v|||+C\kappa)^{p-1}(|||v|||+C\kappa^{(p-1)j_*+1}) g(x) 
\end{align*}
for $x \in {\bf R}^{N}$. 
From $I_{1} + I_{2} = F_{\kappa}[v] \in C_{0}({\bf R}^{N})$, it follows that  
\begin{align} \label{eq:4.13}
|||F_\kappa[v]|||\le C(|||v|||+C\kappa)^{p-1}(|||v|||+C\kappa^{(p-1)j_*+1}) \quad \text{for} \quad 0 < \kappa < 1. 
\end{align}
Since 
\begin{align*}
| F_{\kappa}[v] - F_{\kappa}[\tilde{v}] | = G*|(v + U^{\kappa}_{j_{*}})^{p}_{+} - (\tilde{v} + U^{\kappa}_{j_{*}})^{p}_{+}| 
 \le p G*[(|v| + |\tilde{v}| + U^{\kappa}_{j_{*}})^{p-1} |v - \tilde{v} | ], 
\end{align*}
similarly to \eqref{eq:4.13}, we obtain
$$
|||F_\kappa[v]-F_\kappa[\tilde{v}]|||\le C(|||v|||+|||\tilde{v}|||+C\kappa)^{p-1}|||v-\tilde{v}|||,
\quad v,\tilde{v}\in gBC({\bf R}^N),
$$
for $0<\kappa<1$. 
Then we can find positive constants $\delta$ and $\kappa^0$ such that 
\begin{equation*}
\begin{split}
 & \overline{B_{gBC}(0,\delta)}\ni v\to F_\kappa[v]\in \overline{B_{gBC}(0,\delta)},\\
 & |||F_\kappa[v]-F_\kappa[\tilde{v}]|||\le\frac{1}{2}|||v-\tilde{v}|||
 \quad \text{for} \quad v,\tilde{v}\in\overline{B_{gBC}(0,\delta)},
\end{split}
\end{equation*}
for $0<\kappa<\kappa^0$, where $\overline{B_{gBC}(0,\delta)}:=\{f\in gBC({\bf R}^N)\,:\,|||f|||\le \delta\}$. 

Define a sequence $\{w_k\}\subset gBC({\bf R}^N)$ by 
\begin{equation*}
w_0:=0,\qquad
w_k:=F_\kappa[w_{k-1}]\quad(k=1,2,\dots). 
\end{equation*}
Due to Lemma~\ref{Lemma:3.1}, by induction we see that 
\begin{equation*}
0<w_1<w_2<\cdots<w_k < \cdots \quad \text{in} \quad {\bf R}^{N}. 
\end{equation*}
Applying the fixed point theorem, we can find a positive function $w_\infty\in gBC({\bf R}^N)$ such that 
$$
\lim_{k\to\infty}|||w_\infty-w_k|||=0,
\qquad
w_\infty=F_\kappa[w_\infty]. 
$$
Then, by Lemma~\ref{Lemma:4.1} we see that 
$u=U_{j_*}^\kappa+w_\infty$ is a $(C_0+L^q_{\rm c})$-solution. 
Therefore assertion~(i) follows. 
Assertion~(ii) is proved by Lemma~\ref{Lemma:4.3}. 
The proof is complete.
$\Box$
\vspace{5pt}

Let $u^\kappa$ be a minimal $(C_0+L^q_{\rm c})$-solution of \eqref{eq:1.1}. 
Consider the linearized eigenvalue problem of \eqref{eq:1.1} at $u^\kappa$
\begin{equation} \label{eigen-k}
\tag{$\mbox{E}_\kappa$}
-\Delta\phi+\phi=\lambda p(u^\kappa)^{p-1}\phi\quad\mbox{in}\quad{\bf R}^N,
\qquad
\phi\in H^1({\bf R}^N).
\end{equation}
Since $q>N(p-1)/2$ and $u^\kappa$ is a $(C_0+L^q_{\rm c})$-solution of \eqref{eq:1.1}, 
applying Lemma~\ref{Lemma:2.1}, we see that problem~$(\mbox{E}_\kappa)$ has the first eigenvalue $\lambda_1^\kappa>0$ and 
\begin{equation}
\label{eq:4.14}
\begin{split}
 & \lambda_1^\kappa=\inf\left\{
\|\psi\|_{H^1({\bf R}^N)}^2\,\biggr/\,p\int_{{\bf R}^N}(u^\kappa)^{p-1}\psi^2\,dx\,\right.\\
 & \qquad\qquad\qquad\qquad
\left.:\,
\psi\in H^1({\bf R}^N),\,\,\int_{{\bf R}^N}(u^\kappa)^{p-1}\psi^2\,dx\not=0\right\}.
\end{split}
\end{equation}
Then the following two lemmas hold. 
\begin{lemma}
\label{Lemma:4.5}
Let $\phi^\kappa_1$ be the first eigenfunction of problem $(\mbox{E}_\kappa)$ such that $\phi_1^\kappa>0$ in ${\bf R}^N$. 
Then $\phi_1^\kappa\in C_0({\bf R}^N)$ and 
\begin{equation}
\label{eq:4.15}
\phi_1^\kappa(x)= p \lambda_1^\kappa \int_{{\bf R}^N} G(x-y) (u^\kappa(y))^{p-1}\phi_1^\kappa(y)\,dy,
\qquad x\in{\bf R}^N. 
\end{equation}
Furthermore, there exists $c>0$ such that 
\begin{equation}
\label{eq:4.16}
c^{-1}g(x)\le\phi_1^\kappa(x)\le c g(x)\quad\mbox{in}\quad{\bf R}^N.
\end{equation}
\end{lemma}
{\bf Proof.}
Since $(u^\kappa)^{p-1}\in L^{q/(p-1)}({\bf R}^N)$ with $q/(p-1)>N/2$ and $\phi_1^\kappa\in H^1({\bf R}^N)$, 
we apply regularity theorems for elliptic equations to see that $\phi_1^\kappa\in C_0({\bf R}^N)$. 
Indeed, since $\phi_1^\kappa\in H^1({\bf R}^N)$, it follows that 
$$
\lim_{R\to\infty}\sup_{x\in{\bf R}^N\setminus B(0,R)}\int_{B(x,1)}|\phi_1^\kappa|^2\,dx
\le\lim_{R\to\infty}\int_{{\bf R}^N\setminus B(0,R-1)}|\phi_1^\kappa|^2\,dx=0. 
$$
Then, applying \cite[Theorems~1 and 8]{Serrin} to $(\mbox{E}_\kappa)$, 
we see that $\phi_1^\kappa\in C_0({\bf R}^N)$. 

Set 
$$
\tilde{\phi}(x):= p \lambda_1^\kappa \int_{{\bf R}^N}G(x-y)(u^\kappa(y))^{p-1}\phi_1^\kappa(y)\,dy. 
$$
Then, by a similar argument as in the proof of Lemma~\ref{Lemma:4.2} we see that 
\begin{eqnarray}
\label{eq:4.17}
 & & \mbox{$\tilde{\phi}(x)\le Cg(x)$ in ${\bf R}^N$};\\
\nonumber
 & & \tilde{\phi}\in C_0({\bf R}^N) \cap H^1({\bf R}^N);\\
\nonumber
 & & \mbox{$\tilde{\phi}$ is a weak solution of 
  $-\Delta\phi+\phi= p \lambda_1^\kappa (u^\kappa)^{p-1}\phi_1^\kappa$ in ${\bf R}^N$}. 
\end{eqnarray}
Then we deduce that $z:=\phi_1^\kappa-\tilde{\phi}\in H^1({\bf R}^N)$ and $z$ is a weak solution of 
$-\Delta z+z=0$ in ${\bf R}^N$.
Then it easily follows that $z=0$ in ${\bf R}^N$, i.e., $\phi_1^\kappa=\tilde{\phi}$ in ${\bf R}^N$. 
Thus we have \eqref{eq:4.15}. 

Since $\phi_1^\kappa\in C_0({\bf R}^N)$ and $\phi_1^\kappa>0$ in ${\bf R}^N$, 
Lemma~\ref{Lemma:3.1}~(i) together with \eqref{eq:4.15} implies that 
$$
\phi_1^\kappa(x)\ge C\int_{B(0,1)}G(x-y)\,dy=Cg(x)\quad\mbox{in}\quad{\bf R}^N.
$$
This together with \eqref{eq:4.17} implies \eqref{eq:4.16}. 
Thus Lemma~\ref{Lemma:4.5} follows.
$\Box$
\begin{lemma}
\label{Lemma:4.6}
Define
$$
\kappa^*:=\sup\{\kappa>0\,:\,S_\kappa\not=\emptyset\}. 
$$
Let $0<\kappa<\kappa^{*}$ and let $\lambda^\kappa_1$ be the first eigenvalue to problem~\eqref{eigen-k}. 
Then $\lambda^\kappa_1>1$ and 
\begin{equation}
\label{eq:4.18}
\int_{{\bf R}^N}[|\nabla\psi|^2+\psi^2]\,dx\ge p\lambda^\kappa_1\int_{{\bf R}^N}(u^\kappa)^{p-1}\psi^2\,dx
>p\int_{{\bf R}^N}(u^\kappa)^{p-1}\psi^2\,dx
\end{equation}
for $\psi\in H^1({\bf R}^N)\setminus\{0\}$. 
Furthermore, 
\begin{equation}
\label{eq:4.19}
1<\lambda_1^{\kappa'}\le\lambda_1^\kappa\quad\mbox{if}\quad 0<\kappa\le\kappa'<\kappa^{*}. 
\end{equation}
\end{lemma}
{\bf Proof.}
Let $0<\kappa<\kappa'<\kappa^*$. 
Set 
$$
w^\kappa:=u^\kappa-U_{j_*}^\kappa,
\qquad
w^{\kappa'}:=u^{\kappa'}-U_{j_*}^{\kappa'},
\qquad 
\psi:=w^{\kappa'}-w^\kappa.
$$ 
By Lemma~\ref{Lemma:3.1} we see that 
\begin{equation}
\label{eq:4.20}
\begin{split}
w^{\kappa'} & =\lim_{j\to\infty}[U_j^{\kappa'}-U_{j_*}^{\kappa'}]
=\lim_{j\to\infty}\left[V_j^{\kappa'}+V_{j-1}^{\kappa'}+\dots+V_{j_*+1}^{\kappa'}\right]\\
 & >\lim_{j\to\infty}\left[V_j^\kappa+V_{j-1}^\kappa+\dots+V_{j_*+1}^\kappa\right]
=\lim_{j\to\infty}[U_j^\kappa-U_{j_*}^\kappa]=w^\kappa.
\end{split}
\end{equation}
Since the function 
$$
[0,\infty)\ni s\to (t+s)^p-s^p
$$ 
is strictly monotone increasing for any fixed $t>0$, 
we observe from Lemma~\ref{Lemma:3.1}~(i) that
\begin{equation}
\label{eq:4.21}
\begin{split}
 & (w^{\kappa'}+U_{j_*}^{\kappa'})^p-(U_{j_*-1}^{\kappa'})^p
=(w^{\kappa'}+V_{j_*}^{\kappa'}+U_{j_*-1}^{\kappa'})^p-(U_{j_*-1}^{\kappa'})^p\\
 & \qquad
>(w^{\kappa'}+V_{j_*}^\kappa+U_{j_*-1}^\kappa)^p-(U_{j_*-1}^\kappa)^p
=(w^{\kappa'}+U_{j_*}^\kappa)^p-(U_{j_*-1}^\kappa)^p\ge 0.
\end{split}
\end{equation}
Since $\phi_1^\kappa$ is the eigenfunction of \eqref{eigen-k} and positive,
by \eqref{eq:4.20} we have 
\begin{equation}
\label{eq:4.22}
p\lambda_1^\kappa\int_{{\bf R}^N}(u^\kappa)^{p-1}\psi\phi_1^\kappa\,dx>0. 
\end{equation}
Furthermore, by \eqref{eq:4.4} and \eqref{eq:4.21} we obtain 
\begin{equation}
\label{eq:4.23}
\begin{split}
 & p\lambda_1^\kappa\int_{{\bf R}^N}(u^\kappa)^{p-1}\psi\phi_1^\kappa\,dx
=\int_{{\bf R}^N}[\nabla\psi\nabla\phi_1^\kappa+\psi\phi_1^\kappa]\,dx\\
 & =\int_{{\bf R}^N}[(w^{\kappa'}+U_{j_*}^{\kappa'})^p-(U_{j_*-1}^{\kappa'})^p
-(w^\kappa+U_{j_*}^{\kappa})^p+(U_{j_*-1}^\kappa)^p]\phi_1^\kappa\,dx\\
 & >\int_{{\bf R}^N}[(w^{\kappa'}+U_{j_*}^\kappa)^p
-(w^\kappa+U_{j_*}^{\kappa})^p]\phi_1^\kappa\,dx
\ge p\int_{{\bf R}^N}(u^\kappa)^{p-1}\psi\phi_1^\kappa\,dx.
\end{split}
\end{equation}
Here we used Lemma \ref{Lemma:4.2}. 
Combining \eqref{eq:4.22} with \eqref{eq:4.23},  
we see that $\lambda^\kappa_1>1$. 
This together with \eqref{eq:4.14} implies \eqref{eq:4.18} and \eqref{eq:4.19}. 
Thus Lemma~\ref{Lemma:4.6} follows.
$\Box$
\vspace{5pt}

Now we are ready to complete the proof of Theorem~\ref{Theorem:1.1}. 
\vspace{5pt}
\newline
{\bf Proof of Theorem~\ref{Theorem:1.1}.}
We show that $\kappa^*<\infty$. 
Consider the eigenvalue problem
\begin{equation}
\label{eq:4.24}
-\Delta \psi+\psi=\lambda\psi\quad\mbox{in}\quad B(0,1),
\quad
\psi=0\quad\mbox{on}\quad\partial B(0,1),
\quad
\psi\in H_0^1(B(0,1)). 
\end{equation}
Let $\lambda_B$ and $\psi_B$ 
be the first eigenvalue and the first eigenfunction of problem~\eqref{eq:4.24}, respectively, 
such that $\psi_B>0$ in $B(0,1)$. 
Set $\psi_B=0$ outside $B(0,1)$.  
Then $\psi_B\in H^1({\bf R}^N)$ and 
\begin{equation}
\label{eq:4.25}
\int_{{\bf R}^N}[|\nabla\psi_B|^2+\psi_B^2]\,dx
=\int_{B(0,1)}[|\nabla\psi_B|^2+\psi_B^2]\,dx=\lambda_B\int_{B(0,1)}\psi_B^2\,dx>0.
\end{equation}
On the other hand, it follows from Lemmas~\ref{Lemma:3.1} and \ref{Lemma:4.6} that 
\begin{equation}
\label{eq:4.26}
\begin{split}
\int_{{\bf R}^N}[|\nabla\psi_B|^2+\psi_B^2]\,dx
 & >p\int_{{\bf R}^N}(u^\kappa)^{p-1}\psi_B^2\,dx\\
 & \ge p\int_{{\bf R}^N}(U_0^\kappa)^{p-1}\psi_B^2\,dx
=p\kappa^{p-1}\int_{{\bf R}^N}\mu_0^{p-1}\psi_B^2\,dx
\end{split}
\end{equation}
for $0<\kappa<\kappa^{*}$. 
By \eqref{eq:3.1}, \eqref{eq:4.25} and \eqref{eq:4.26} we see that 
$$
\kappa^{p-1}\le\frac{\lambda_B}{p}
\int_{B(0,1)}\psi_B^2\,dx\biggr/\int_{{\bf R}^N}\mu_0^{p-1}\psi_B^2\,dx
$$
for $0<\kappa<\kappa^*$. Then we deduce that $\kappa^*<\infty$. 
This together with Lemma~\ref{Lemma:4.4} implies 
assertions~(i) and (ii) of Theorem~\ref{Theorem:1.1}. 
Thus the proof of Theorem~\ref{Theorem:1.1} is complete.
$\Box$
%%%%%%%%%%%%%%%%%%%%%%%%%%%%%%%%%%%%%
%%%%%%%%%%%%%%%%%%%%%%%%%%%%%%%%%%%%%
\section{Uniform estimates of $w^\kappa$} \label{section:5}
%%%%%%%%%%%%%%%%%%%%%%%%%%%%%%%%%%%%%
%%%%%%%%%%%%%%%%%%%%%%%%%%%%%%%%%%%%%
For $0<\kappa<\kappa^*$, 
let $u^\kappa$ be the minimal $(C_0+L^q_{\rm c})$-solution of \eqref{eq:1.1} 
and set $w^\kappa:=u^\kappa-U_{j_*}^\kappa$. 
In this section we obtain uniform estimates of $\{w^\kappa\}_{0<\kappa<\kappa^*}$ 
in $H^1({\bf R}^N)$ and $L^\infty({\bf R}^N)$ and prove the following proposition:
\begin{proposition}
\label{Proposition:5.1}
Let $1<p<p_{JL}$ and 
assume the same conditions as in Theorem~{\rm\ref{Theorem:1.1}}. 
Then problem~\eqref{eq:1.1} with $\kappa=\kappa^*$ possesses a $(C_0+L^q_{\rm c})$-solution. 
\end{proposition}
We prepare an inequality in order to estimate the nonlinear term of \eqref{eq:4.3}. 
\begin{lemma}
\label{Lemma:5.1}
Let $p>1$, $\epsilon>0$ and $\delta\ge 0$. 
Then there exists a constant $c > 0$ such that 
\begin{equation}
\label{eq:5.1}
t^p-s^p
\le (1+\epsilon)t^{p-1}(t-s)+ c(t-s)^{1-\delta}s^{p-1+\delta} \quad \text{for} \quad t\ge s\ge 0. 
\end{equation}
\end{lemma}
{\bf Proof.}
Let $\delta>0$ and $\epsilon>0$. 
For sufficiently small $\beta > 0$, we have
\begin{equation*}
\begin{split}
 s[t^{p-1}-s^{p-1}]
 & =(p-1)s\int_s^t \tau^{(p-1)(1-\beta)}\tau^{-\beta\delta}\tau^{(p-1)\beta-1+\beta\delta}\,d\tau\\
 & \le Ct^{(p-1)(1-\beta)}[t^{1-\beta\delta}-s^{1-\beta\delta}]s^{(p-1)\beta+\beta\delta}\\
 & \le Ct^{(p-1)(1-\beta)}(t-s)^{1-\beta\delta}s^{(p-1)\beta+\beta\delta}\\
 & \le (t-s)\left[\epsilon t^{p-1}+C(t-s)^{-\delta}s^{p-1+\delta}\right]
\end{split}
\end{equation*}
for $t>s\ge 0$. This implies that 
\begin{equation*}
\begin{split}
 t^p-s^p
 & =t^{p-1}(t-s)+t^{p-1}s-s^p\\
 & =t^{p-1}(t-s)+s[t^{p-1}-s^{p-1}]\\
 & \le (1+\epsilon)t^{p-1}(t-s)+C(t-s)^{1-\delta}s^{p-1+\delta}
\end{split}
\end{equation*}
for $t\ge s\ge 0$. Thus Lemma~\ref{Lemma:5.1} follows.
$\Box$
\vspace{5pt}

We obtain a uniform estimate of $\{w^\kappa\}_{0<\kappa<\kappa^*}$ in $H^1({\bf R}^N)$.
\begin{lemma}
\label{Lemma:5.2}
Let $0<\kappa<\kappa^{*}$. 
Let $u^\kappa$ be the minimal $(C_0+L^{q}_{\rm c})$-solution of \eqref{eq:1.1} 
and set $w^\kappa:=u^\kappa-U_{j_*}^\kappa$. 
Then 
\begin{equation}
\label{eq:5.2}
\sup_{0<\kappa<\kappa^*}\|w^\kappa\|_{H^1({\bf R}^N)}<\infty.  
\end{equation}
\end{lemma}
{\bf Proof.}
Let $\epsilon$ and $\delta$ be sufficiently small positive constants. 
By \eqref{eq:3.10} we can find $\tilde{V}\in H^1({\bf R}^N)$ such that $V_{j_*}^{\kappa^*}\le\tilde{V}$ in ${\bf R}^N$. 
Then it follows from \eqref{eq:4.4} and \eqref{eq:5.1} that 
\begin{equation*}
\begin{split}
\|w^\kappa\|_{H^1({\bf R}^N)}^2 & =\int_{{\bf R}^N}[(w^\kappa+U_{j_*}^\kappa)^p-(U_{j_*-1}^\kappa)^p]w_\kappa\,dx\\
 & \le(1+\epsilon)\int_{{\bf R}^N}(w^\kappa+U_{j_*}^\kappa)^{p-1}(w_\kappa+V_{j_*}^\kappa)w^\kappa\,dx\\
 & \qquad\quad
 +C\int_{{\bf R}^N}(U_{j_*-1}^\kappa)^{p-1+\delta}(w^\kappa + V_{j_{*}}^{\kappa})^{1-\delta}w^\kappa\,dx\\
 & \le(1+\epsilon)^2\int_{{\bf R}^N}(u^\kappa)^{p-1}[(w^\kappa)^2+C\tilde{V}^2]\,dx\\
 & \qquad\quad
 +C\int_{{\bf R}^N}(U_{j_*}^{\kappa^*})^{p-1+\delta}\left((w^\kappa)^{2-\delta}+ M^{2-\delta} \right)\,dx,
\end{split}
\end{equation*}
where $M:=\|V_{j_*}^{\kappa^*}\|_{L^\infty({\bf R}^N)}$. 
This together with Lemma~\ref{Lemma:4.6} implies that 
\begin{equation}
\label{eq:5.3}
\begin{split}
\|w^\kappa\|_{H^1({\bf R}^N)}^2 & \le\frac{(1+\epsilon)^2}{p}
\left[\|w_\kappa\|_{H^1({\bf R}^N)}^2+C\|\tilde{V}\|_{H^1({\bf R}^N)}^2\right]\\
 & \qquad\qquad
+C\int_{{\bf R}^N}(U_{j_*}^{\kappa^*})^{p-1+\delta}(w_\kappa)^{2-\delta}\,dx+C.
\end{split}
\end{equation}
On the other hand, 
it follows from $q>\max\{N(p-1)/2,p\}>1$ that 
\begin{equation}
\label{eq:5.4}
2<\frac{2q}{q-p+1}<p_S+1.
\end{equation}
By \eqref{eq:3.10}, 
applying the interpolation inequality in Lebesgue spaces 
and the Sobolev inequality, we obtain
\begin{equation}
\label{eq:5.5}
\begin{split}
 & \int_{{\bf R}^N}(U_{j_*}^\kappa)^{p-1+\delta}(w_\kappa)^{2-\delta}\,dx
\le\|U_{j_*}^{\kappa^*}\|_{L^q({\bf R}^N)}^{p-1+\delta}\|w^\kappa\|^{2-\delta}_{L^{\frac{(2-\delta)q}{q-p+1-\delta}}({\bf R}^N)}\\
 & \le C\|U_{j_*}^{\kappa^*}\|_{L^q({\bf R}^N)}^{p-1+\delta}[\|w^\kappa\|_{L^{p_S+1}({\bf R}^N)}+\|w^\kappa\|_2]^{2-\delta}\\
 & \le C\|U_{j_*}^{\kappa^*}\|_{L^q({\bf R}^N)}^{p-1+\delta}\|w^\kappa\|^{2-\delta}_{H^1({\bf R}^N)}
\le\epsilon\|w^\kappa\|_{H^1({\bf R}^N)}^2+C.
\end{split}
\end{equation}
Taking a sufficiently small $\epsilon>0$ if necessary, 
by \eqref{eq:5.3} and \eqref{eq:5.5} we obtain 
$$
\|w^\kappa\|_{H^1({\bf R}^N)}^2\le\frac{1}{2} \left(1+\dfrac{1}{p} \right)\|w^\kappa\|_{H^1({\bf R}^N)}^2+C.
$$
This implies \eqref{eq:5.2}. Thus Lemma~\ref{Lemma:5.2} follows.
$\Box$
\vspace{5pt}

Applying Lemma~\ref{Lemma:5.2}, we have:
\begin{lemma}
\label{Lemma:5.3}
Assume the same conditions as in Theorem~{\rm \ref{Theorem:1.1}}.
Let $0<\kappa<\kappa^{*}$ and $z\in{\bf R}^N$. 
Let $\nu\in(0,1)$ be such that 
\begin{equation}
\label{eq:5.6}
4\nu(1-\nu)p>1. 
\end{equation}
Then 
\begin{equation}
\label{eq:5.7}
\sup_{z\in{\bf R}^N}\sup_{0<\kappa<\kappa^*}\|w^\kappa g_z\|_{L^{\frac{p_S+1}{2\nu}}({\bf R}^N)}<\infty,
\end{equation}
where $g_z(x):=g(x-z)$ for $x,z\in{\bf R}^N$. 
\end{lemma}
{\bf Proof.}
Let $0<\kappa<\kappa^{*}$, $z\in{\bf R}^N$ and $0<\delta<1$. 
Similarly to the proof of Lemma~\ref{Lemma:5.2}, 
we set $M:=\|V_{j_*}^{\kappa^*}\|_{L^\infty({\bf R}^N)}$. 
Then it follows from Lemma~\ref{Lemma:3.1}~(ii) that 
\begin{equation*}
M\ge\sup_{0<\kappa\le\kappa^*}\|V_{j_*}^\kappa\|_{L^\infty({\bf R}^N)}. 
\end{equation*}
Set 
$$
\Psi(x):=\left[(M+w^\kappa(x))g_z(x)\right]^{\frac{1}{2\nu}}. 
$$
For any $\epsilon>0$, 
it follows that 
\begin{equation*}
\begin{split}
 & \int_{{\bf R}^N}|\nabla\Psi|^2\,dx
 =\int_{{\bf R}^N}\left|\frac{1}{2\nu}(M+w^\kappa)^{\frac{1}{2\nu}-1}
 g_z^{\frac{1}{2\nu}}[\nabla w^\kappa +(M+w^\kappa)g_z^{-1}\nabla g_z]\right|^2\,dx\\
 & \le\int_{{\bf R}^N}\left(\frac{1}{4\nu^2}+\epsilon\right)
 (M+w^\kappa)^{\frac{1}{\nu}-2}g_z^{\frac{1}{\nu}}|\nabla w^\kappa|^2\,dx
 +C\nu^{-2}\int_{{\bf R}^N}(M+w^\kappa)^{\frac{1}{\nu}}g_z^{\frac{1}{\nu}-2}|\nabla g_z|^2\,dx
\end{split}
\end{equation*}
and
\begin{equation*}
\begin{split}
 & \int_{{\bf R}^N}\nabla w^\kappa\nabla \left[(M+w^\kappa)^{1-\nu}g_z\right]^{\frac{1}{\nu}}\,dx\\
 & =\frac{1}{\nu}\int_{{\bf R}^N}(M+w^\kappa)^{\frac{1}{\nu}-2}g_z^{\frac{1}{\nu}}
 [(1-\nu)\nabla w^\kappa+(M+ w^{\kappa})g_z^{-1}\nabla g_z]\nabla w^\kappa\,dx\\
 & \ge\left(\frac{1-\nu}{\nu}-\epsilon\right)\int_{{\bf R}^N}(M+w^\kappa)^{\frac{1}{\nu}-2}g_z^{\frac{1}{\nu}}|\nabla w^\kappa|^2\,dx
 -C\nu^{-1}\int_{{\bf R}^N}(M+w^\kappa)^{\frac{1}{\nu}}g_z^{\frac{1}{\nu}-2}|\nabla g_z|^2\,dx.
\end{split}
\end{equation*}
Then, taking $\epsilon>0$ small enough, we observe from \eqref{eq:2.2} that 
\begin{equation}
\label{eq:5.8}
4(\nu-\nu^2-\delta)\|\Psi\|_{H^1({\bf R}^N)}^2
\le\left(w^\kappa,\left[(M+w^\kappa)^{1-\nu}g_z\right]^{\frac{1}{\nu}}\right)_{H^1({\bf R}^N)}
+C\|\Psi\|_2^2.
\end{equation}
On the other hand, 
it follows from \eqref{eq:3.3}, \eqref{eq:4.4} and \eqref{eq:5.1} that 
\begin{equation}
\label{eq:5.9}
\begin{split}
 & \left(w^\kappa,\left[(M+w^\kappa)^{1-\nu}g_z\right]^{\frac{1}{\nu}}\right)_{H^1({\bf R}^N)}\\
 & =\int_{{\bf R}^N}[(w^{\kappa}+U_{j_*}^\kappa)^p-(U_{j_*-1}^\kappa)^p]\left[(M+w^\kappa)^{1-\nu}g_z\right]^{\frac{1}{\nu}}\,dx\\
 & \le(1+\delta)\int_{{\bf R}^N}(w^{\kappa}+U_{j_*}^\kappa)^{p-1}(w^{\kappa}+V_{j_*}^\kappa)
 \left[(M+w^\kappa)^{1-\nu}g_z\right]^{\frac{1}{\nu}}\,dx\\
 & \qquad\quad
 +C\int_{{\bf R}^N}(U_{j_*-1}^\kappa)^{p-1}
 (w^{\kappa}+V_{j_*}^\kappa)\left[(M+w^\kappa)^{1-\nu}g_z\right]^{\frac{1}{\nu}}\,dx\\
  & \le(1+\delta)\int_{{\bf R}^N}(u^\kappa)^{p-1}\Psi^2\,dx+C\int_{{\bf R}^N}(U_{j_{*}}^{\kappa^{*}})^{p-1}\Psi^2\,dx.
\end{split}
\end{equation}
Furthermore, Lemma~\ref{Lemma:4.6} implies that
\begin{equation}
\label{eq:5.10}
\begin{split}
\int_{{\bf R}^N}(u^\kappa)^{p-1}\Psi^2\,dx<\frac{1}{p}\|\Psi\|_{H^1({\bf R}^N)}^2.
\end{split}
\end{equation}
By the H\"older inequality and Lemma~\ref{Lemma:3.2}
we have
\begin{equation}
\label{eq:5.11}
\begin{split}
\int_{{\bf R}^N}(U_{j_*}^{\kappa^{*}})^{p-1}\Psi^2\,dx
 \le \|U_{j_*}^{\kappa^{*}}\|_{L^q({\bf R}^N)}^{p-1} \|\Psi\|^2_{L^{\frac{2q}{q-p+1}}({\bf R}^N)}
 \le C\|\Psi\|^2_{L^{\frac{2q}{q-p+1}}({\bf R}^N)}.
\end{split}
\end{equation}
Therefore, by \eqref{eq:5.8}, \eqref{eq:5.9}, \eqref{eq:5.10} and \eqref{eq:5.11} 
we obtain
$$
4(\nu-\nu^2-\delta)\|\Psi\|_{H^1({\bf R}^N)}^2
\le\frac{1+\delta}{p}\|\Psi\|_{H^1({\bf R}^N)}^2+C\|\Psi\|_{L^2({\bf R}^N)}^2+C\|\Psi\|^2_{L^{\frac{2q}{q-p+1}}({\bf R}^N)}. 
$$
By \eqref{eq:5.6}, taking a sufficiently small $\delta>0$ if necessary, 
we see that 
$$
\|\Psi\|_{H^1({\bf R}^N)}^2\le C\|\Psi\|_{L^2({\bf R}^N)}^2+C\|\Psi\|^2_{L^{\frac{2q}{q-p+1}}({\bf R}^N)}.
$$
Then we deduce from the Sobolev inequality that 
$$
\|\Psi\|_{L^{p_S+1}({\bf R}^N)}^2\le C\|\Psi\|_{L^2({\bf R}^N)}^2+C\|\Psi\|^2_{L^{\frac{2q}{q-p+1}}({\bf R}^N)}.
$$
This implies that 
\begin{equation}
\label{eq:5.12}
\|w^\kappa g_z\|_{L^{\frac{p_S+1}{2\nu}}({\bf R}^N)}
\le C\|w^\kappa g_z\|_{L^{\frac{1}{\nu}}({\bf R}^N)}+C\|w^\kappa g_z\|_{L^{\frac{q}{\nu(q-p+1)}}({\bf R}^N)} + C. 
\end{equation}

Consider the case where $1/\nu\ge 2$. 
By \eqref{eq:5.4}, applying the interpolation inequality in Lebesgue spaces, 
for any $\epsilon'>0$, we can find a constant $C>0$ such that
\begin{equation}
\label{eq:5.13}
\|w^\kappa g_z\|_{L^{\frac{1}{\nu}}({\bf R}^N)}+\|w^\kappa g_z\|_{L^{\frac{q}{\nu(q-p+1)}}({\bf R}^N)}
\le\epsilon'\|w^\kappa g_z\|_{L^{\frac{p_S+1}{2\nu}}({\bf R}^N)}+C\|w^\kappa g_z\|_{L^2({\bf R}^N)}.
\end{equation}
This together with \eqref{eq:5.12} implies that 
\begin{equation}
\label{eq:5.14}
\|w^\kappa g_z\|_{L^{\frac{p_S+1}{2\nu}}({\bf R}^N)} 
 \le C\|w^\kappa g_z\|_{L^2({\bf R}^N)} + C 
 \le C\|w^\kappa\|_{L^2({\bf R}^N)} + C. 
\end{equation}
In the case where $1/\nu<2$, similarly to \eqref{eq:5.13}, we have 
$$
\|w^\kappa g_z\|_{L^{\frac{q}{\nu(q-p+1)}}({\bf R}^N)}
\le\epsilon'\|w^\kappa g_z\|_{L^{\frac{p_S+1}{2\nu}}({\bf R}^N)}+C\|w^\kappa g_z\|_{L^{\frac{1}{\nu}}({\bf R}^N)}.
$$
Then, by \eqref{eq:5.12} we have
\begin{equation}
\label{eq:5.15}
\begin{split}
\|w^\kappa g_z\|_{L^{\frac{p_S+1}{2\nu}}({\bf R}^N)} 
 & \le C \|w^\kappa g_z\|_{L^{\frac{1}{\nu}}({\bf R}^N)} + C 
\le C \left\|[(w^\kappa)^{2\nu}+1]g_z\right\|_{L^{\frac{1}{\nu}}({\bf R}^N)} + C \\
 & \le C\|w^\kappa\|_{L^2({\bf R}^N)}^{2/\nu}+C. 
\end{split}
\end{equation}
Therefore, 
by Lemma~\ref{Lemma:5.2}, \eqref{eq:5.14} and \eqref{eq:5.15}, 
we obtain \eqref{eq:5.7}. Thus Lemma~\ref{Lemma:5.3} follows. 
$\Box$\vspace{5pt}

We apply regularity theorem for elliptic equations 
to obtain the following lemma.
\begin{lemma}
\label{Lemma:5.4}
Let $0<\kappa<\kappa^{*}$ and $z\in{\bf R}^N$. 
Assume that there exists $\nu\in(0,1)$ such that 
\begin{equation}
\label{eq:5.16}
4\nu(1-\nu)p>1,
\qquad
\frac{p_S+1}{2\nu}>\frac{N}{2}(p-1).
\end{equation}
Then $\{w^\kappa\}_{0<\kappa<\kappa^*}$ are uniformly bounded in ${\bf R}^N$ and 
equi-continuous for any compact set in ${\bf R}^N$.  
\end{lemma}
{\bf Proof.}
Let $0<\kappa<\kappa^*$. 
It follows that 
\begin{equation*}
\begin{split}
0 & \le (w^\kappa+U_{j_*}^\kappa)^p-(U_{j_*-1}^\kappa)^p
\le p(w^\kappa+U_{j_*}^\kappa)^{p-1}(w^\kappa+M)\\
 & \le C[(w^\kappa)^{p-1}+(U_{j_*}^{\kappa^*})^{p-1}](w^\kappa+M), 
\end{split}
\end{equation*}
where $M:= \| V^{\kappa^{*}}_{j_{*}} \|_{L^{\infty}({\bf R}^{N})}$. 
By Lemma~\ref{Lemma:5.3} and \eqref{eq:5.16} we see that
$$
\sup_{0<\kappa<\kappa^*}\|(w^\kappa)^{p-1}\|_{L^{\frac{p_S+1}{2\nu(p-1)}}({\bf R}^N)}<\infty
\quad\mbox{and}\quad
\frac{p_S+1}{2\nu(p-1)}>\frac{N}{2}.
$$
Furthermore, it follows from Lemma~\ref{Lemma:3.2} that
$$
(U_{j_*}^{\kappa^*})^{p-1}\in L^{\frac{q}{p-1}}({\bf R}^N)
\quad\mbox{and}\quad
\frac{q}{p-1}>\frac{N}{2}. 
$$ 
Then, applying \cite[Theorem~1]{Serrin} to problem~\eqref{eq:4.3},
we obtain 
$$
\|w^\kappa\|_{L^\infty(B(z,1))}
\le C\|w^\kappa\|_{L^2(B(z,2))}+C
$$
for $z\in{\bf R}^N$ and $0<\kappa<\kappa^*$.
This together with Lemma~\ref{Lemma:5.2} implies that 
$\{w^\kappa\}_{0<\kappa<\kappa^*}$ are uniformly bounded in ${\bf R}^N$. 
Furthermore, applying \cite[Theorem~7]{Serrin} to problem~\eqref{eq:4.3}, 
we see that $\{w^\kappa\}_{0<\kappa<\kappa^*}$ are
equi-continuous for any compact set in ${\bf R}^N$. 
Thus Lemma~\ref{Lemma:5.4} follows.
$\Box$
\vspace{5pt}

We show that \eqref{eq:5.16} holds for some $\nu\in(0,1)$ if $1<p<p_{JL}$. 
\begin{lemma}
\label{Lemma:5.5}
Let $1<p<p_{JL}$. 
Then there exists $\nu\in(0,1)$ satisfying \eqref{eq:5.16}. 
\end{lemma}
{\bf Proof.}
Let $1<p<p_{JL}$. 
Let $\nu^\pm$ be the roots of 
$4p(\nu-\nu^2)=1$ such that $\nu^-\le \nu^+$, that is
$$
\nu^{\pm} =\frac{p\pm\sqrt{p^2-p}}{2p}.
$$
Then $4p(\nu-\nu^2)>1$ if and only if $\nu^{-}<\nu<\nu^{+}$. 
Since $\nu^-<(0,1/2)$, 
if 
\begin{equation}
\label{eq:5.17}
\frac{p_S+1}{2\nu^{-}}>\frac{N}{2}(p-1),
\end{equation}
then we can find $\nu\in(0,1)$ satisfying \eqref{eq:5.16}. 
 
We prove \eqref{eq:5.17}.
If $N=2$, then $p_S=\infty$ and \eqref{eq:5.17} holds. So it suffices to consider the case of $N\ge 3$. 
Since
$$
\frac{N}{\nu^{-}(N-2)}=\frac{p_{S}+1}{2\nu^{-}}>\frac{N}{2}(p-1),
$$
the inequality \eqref{eq:5.17} is equivalent to 
%it follows that 
$$
\frac{2}{N-2}>\nu^{-}(p-1)=\frac{p-1}{2}-\frac{p-1}{2p}\sqrt{p^2-p}. 
$$
This is reduced to 
%This implies that
\begin{equation}
\label{eq:5.18}
\sqrt{p^2-p}>\frac{2p}{p-1}\left[\frac{p-1}{2}-\frac{2}{N-2}\right]
=p-\frac{4p}{(p-1)(N-2)}=:D.
\end{equation}
If $p<p_S$, then $D<0$ and we have \eqref{eq:5.18}. 
Thus \eqref{eq:5.17} holds in the case of $p<p_S$. 

Consider the case of $p\ge p_S$. 
Then $D\ge 0$ and \eqref{eq:5.18} is equivalent to 
$$
p^2-p>D^2=p^2 -\frac{8p^2}{(p-1)(N-2)}+\frac{16p^2}{(p-1)^2(N-2)^2},
$$
that is  
\begin{equation}
\label{eq:5.19}
(N-2)(N-10)p^2-2(N^2-8N+4)p+(N-2)^2<0.
\end{equation}

Let $N=10$. 
Then \eqref{eq:5.19} is equivalent to $p>4/3$. 
Since $p\ge p_S=(N+2)/(N-2)=3/2$, \eqref{eq:5.19} holds for $p\ge p_S$ if $N=10$. 

Let $N\not=10$. 
Let $p^\pm$ be the roof of 
$$
(N-2)(N-10)p^2-2(N^2-8N+4)p+(N-2)^2=0
$$
such that $p^-\le p^+$, that is 
$$
p^\pm=\frac{N^2-8N+4\pm8\sqrt{N-1}}{(N-2)(N-10)}.
$$
In the case of $N\ge 11$, 
since $p_S\le p<p_{JL}=p^+$, 
\eqref{eq:5.19} holds if $p^- < p_S$. 
Furthermore, $p^- < p_S$ if and only if 
$$
4-8\sqrt{N-1}<-20,
$$
which holds if $N\ge 11$. 
Thus \eqref{eq:5.19} holds in the case of  $N\ge 11$. 
In the case of $3\le N\le 9$, 
$(N-2)(N-10)<0$ and $p^+<p_S$. 
Then we have \eqref{eq:5.19} for $p\ge p_S$.  
Therefore \eqref{eq:5.19} holds for $p\ge p_S$ in the case of $N\not=10$. 
Thus Lemma~\ref{Lemma:5.5} follows.
$\Box$
\vspace{5pt}

We are ready to prove Proposition~\ref{Proposition:5.1}.
\vspace{5pt}
\newline
{\bf Proof of Proposition~\ref{Proposition:5.1}.}
Let $1<p<p_{JL}$. Thanks to Lemma~\ref{Lemma:5.5}, we can find $\nu\in(0,1)$ satisfying \eqref{eq:5.16}.  
Therefore, applying Lemmas~\ref{Lemma:5.2} and \ref{Lemma:5.4}, we obtain
$$
\sup_{0<\kappa<\kappa^*}\,\|w^\kappa\|_{H^1}<\infty,
\quad
\sup_{0<\kappa<\kappa^*}\,\|w^\kappa\|_\infty<\infty.
$$
Furthermore, $\{w^\kappa\}_{0<\kappa<\kappa^*}$ are
equi-continuous for any compact set in ${\bf R}^N$. 
Then we can find a solution $w_*\in H^1({\bf R}^N)\,\cap\,C({\bf R}^N)$ of \eqref{eq:4.3} with $\kappa=\kappa^*$. 
Furthermore, similarly to the proof of Lemma~\ref{Lemma:4.5}, 
applying \cite[Theorem~1]{Serrin} to \eqref{eq:4.3}, 
we see that $w_*\in C_0({\bf R}^N)$. 
Therefore, by Lemma~\ref{Lemma:4.1} we can find a $(C_0+L^q_{\rm c})$-solution of \eqref{eq:1.1} 
with $\kappa=\kappa^*$. 
Thus Proposition~\ref{Proposition:5.1} follows. 
$\Box$
%%%%%%%%%%%%%%%%%%%%%%%%%%%%%%%%%%%%%
%%%%%%%%%%%%%%%%%%%%%%%%%%%%%%%%%%%%%
\section{Proof of Theorem~\ref{Theorem:1.2}} \label{section:6}
%%%%%%%%%%%%%%%%%%%%%%%%%%%%%%%%%%%%%
%%%%%%%%%%%%%%%%%%%%%%%%%%%%%%%%%%%%%
Let $u^{\kappa^{*}}$ be the minimal solution of \eqref{eq:1.1} with $\kappa=\kappa^{*}$ 
obtained by Proposition \ref{Proposition:5.1}. 
For the proof of Theorem~\ref{Theorem:1.2}, 
we prepare two lemmas. 
\begin{lemma}
\label{Lemma:6.1}
Let $1<p<p_{JL}$ and assume the same conditions as in Theorem~{\rm\ref{Theorem:1.1}}. 
Then $\lambda_1^{\kappa^*}\ge 1$. 
\end{lemma}
{\bf Proof.}
It follows from \eqref{eq:4.18} that 
$$
\int_{{\bf R}^N}[|\nabla\psi|^2+\psi^2]\,dx
>p\int_{{\bf R}^N}(u^\kappa)^{p-1}\psi^2\,dx
\ge p\int_{{\bf R}^N}(U^\kappa_j)^{p-1}\psi^2\,dx
$$
for $\psi\in H^1({\bf R}^N)\setminus\{0\}$, $\kappa\in(0,\kappa^*)$ and $j\in\{0,1,2,\dots\}$.
Then we have
$$
\int_{{\bf R}^N}[|\nabla\psi|^2+\psi^2]\,dx
\ge p\int_{{\bf R}^N}(U^{\kappa^{*}}_j)^{p-1}\psi^2\,dx
$$
for $\psi\in H^1({\bf R}^N)$ and $j\in\{0,1,2,\dots\}$.
This together with Lemma~\ref{Lemma:3.1}~(ii) implies that 
$$
\int_{{\bf R}^N}[|\nabla\psi|^2+\psi^2]\,dx
\ge p\int_{{\bf R}^N}(u^{\kappa^{*}})^{p-1}\psi^2\,dx
$$
for $\psi\in H^1({\bf R}^N)$. 
Then we see that $\lambda_1^{\kappa^*}\ge 1$, and the proof is complete.
$\Box$
\begin{lemma}
\label{Lemma:6.2}
Let $1<p<p_{JL}$ and assume the same conditions as in Theorem~{\rm\ref{Theorem:1.1}}. 
Let $\kappa\in(0,\kappa^*]$ be such that $\lambda_1^\kappa>1$. 
Then $\kappa<\kappa^*$. 
\end{lemma}
{\bf Proof.}
Let $\kappa\in(0,\kappa^*]$ be such that $\lambda_1^\kappa>1$. 
For $\epsilon\in(0,1]$ and $a\in(0,1]$, set
$$
\tilde{w}^\kappa:=u^\kappa-U_{j_*}^\kappa,
\qquad
\overline{u}:=U_{j_*}^{\kappa+\epsilon}+\tilde{w}^{\kappa}+a\phi_1^\kappa\in C_0({\bf R}^N)+L^q_{\rm c}({\bf R}^N).
$$
Then 
\begin{equation}
\label{eq:6.1}
\begin{split}
 & \overline{u}-G*\overline{u}^p-(\kappa+\epsilon)G*\mu\\
 & =(\tilde{w}^\kappa+a\phi^{\kappa}_{1})-G*\overline{u}^p+G*(U_{j_{*}-1}^{\kappa+\epsilon})^p\\
 & =G*[(\tilde{w}^\kappa+U_{j_*}^\kappa)^p-(U_{j_{*}-1}^\kappa)^p]-G*\overline{u}^p+G*(U_{j_{*}-1}^{\kappa+\epsilon})^p+a\phi_1^\kappa\\
 & =-G*[H(1)-H(0)]+a\phi_1^\kappa,
\end{split}
\end{equation}
where
\begin{equation*}
\begin{split}
 & H(x,t):=P(x,t)^p-Q(x,t)^p,\\
 & P(x,t):=tU_{j_*}^{\kappa+\epsilon}(x)+(1-t)U_{j_*}^{\kappa}(x)+\tilde{w}^{\kappa}(x)+at\phi_1^\kappa(x),\\
 & Q(x,t):=tU_{j_{*}-1}^{\kappa+\epsilon}(x)+(1-t)U_{j_{*}-1}^{\kappa}(x). 
\end{split}
\end{equation*}
It follows that 
\begin{equation}
\label{eq:6.2}
\begin{split}
H'(t) & =pP(x,t)^{p-1}
(U_{j_*}^{\kappa+\epsilon}-U_{j_*}^\kappa+a\phi_1^\kappa)
-pQ(x,t)^{p-1}(U_{j_{*}-1}^{\kappa+\epsilon}-U_{j_{*}-1}^{\kappa})\\
 & =p\left[U_{j_*}^{\kappa}+\tilde{w}^{\kappa}\right]^{p-1}a\phi_1^\kappa\\
 & \quad
 +p\left\{P(t)^{p-1}-\left[U_{j_*}^{\kappa}+\tilde{w}^{\kappa}\right]^{p-1}\right\}a\phi_1^\kappa
 +pP(t)^{p-1}(V_{j_*}^{\kappa+\epsilon}-V_{j_*}^\kappa)\\
 & \quad
 +p\left\{P(t)^{p-1}-Q(t)^{p-1}\right\}(U_{j_{*}-1}^{\kappa+\epsilon}-U_{j_{*}-1}^{\kappa})
\end{split}
\end{equation}
for $x\in{\bf R}^N$ and $t\in(0,1)$.

It follows that
\begin{equation}
\label{eq:6.3}
0<\tilde{w}(x)\le u^\kappa(x)-U_{j_{*}-1}^\kappa(x)=w^\kappa(x)\le Cg(x)\quad\mbox{in}\quad{\bf R}^N.
\end{equation}
Combing Lemmas~\ref{Lemma:3.1}, \ref{Lemma:3.2} and \ref{Lemma:3.3} with \eqref{eq:3.3} and \eqref{eq:6.3}, 
we have
\begin{equation}
\label{eq:6.4}
\begin{split}
 & 0 \le P(x,t)-Q(x,t)
=tV_{j_*}^{\kappa+\epsilon}(x)+(1-t)V_{j_*}^\kappa(x)+\tilde{w}^\kappa(x)+at\phi_1^\kappa(x)\\
 & \quad
 \le V_{j_*}^{\kappa+1}(x)+\tilde{w}^\kappa(x)+\phi_1^\kappa(x)\le Cg(x)\le C\phi_1^\kappa(x),\\
 & \max\{P(x,t),Q(x,t)\}\le U_{j_*}^{\kappa+\epsilon}(x)+\tilde{w}^\kappa(x)+at\phi_1^\kappa(x)\\
 & \quad
 \le CU_{j_*}^\kappa(x)+Cg(x)=C[U_{j_{*}-1}^\kappa(x)+V_{j_*}^\kappa(x)]+Cg(x)\\
 & \quad
 \le CU_{j_{*}-1}^\kappa(x)+Cg(x)\le CU_{j_{*}-1}^\kappa(x),\\
 & \min\{P(x,t),Q(x,t)\}\ge U_{j_{*}-1}^\kappa(x).
\end{split}
\end{equation}
These together with Lemma~\ref{Lemma:3.3} imply that 
\begin{equation}
\label{eq:6.5}
\begin{split}
0 & \le p\left\{P(x,t)^{p-1}-Q(x,t)^{p-1}\right\}(U_{j_{*}-1}^{\kappa+\epsilon}(x)-U_{j_{*}-1}^{\kappa}(x))\\
 & \le C\epsilon(U_{j_{*}-1}^\kappa(x))^{p-1}\phi_1^\kappa(x)
 \le C\epsilon(u^\kappa(x))^{p-1}\phi_1^\kappa(x) 
\end{split}
\end{equation}
for $x\in{\bf R}^N$ and $t\in(0,1)$.
Furthermore, by Lemma~\ref{Lemma:3.4} we have
\begin{equation}
\label{eq:6.6}
\begin{split}
 & pP(x,t)^{p-1}(V_{j_*}^{\kappa+\epsilon}-V_{j_*}^\kappa)\\
 & \le C\epsilon^{\min\{p-1,1\}}(U_{j_{*}-1}^\kappa(x))^{p-1}g(x)
 \le C\epsilon^{\min\{p-1,1\}}(u^\kappa(x))^{p-1}\phi_1^\kappa(x) 
\end{split}
\end{equation}
for $x\in{\bf R}^N$ and $t\in(0,1)$.
Similarly to \eqref{eq:6.4}, we have
\begin{equation*}
\begin{split}
 & 0\le P(x,t)-[U_{j_*}^\kappa(x)+\tilde{w}^\kappa(x)]
=t[U_{j_*}^{\kappa+\epsilon}(x)-U_{j_*}^\kappa(x)]+at\phi_1^\kappa(x)\\
 & \,\,\,\,\le C\epsilon U_{j_*}^\kappa(x)+Cag(x)\le C(\epsilon+a) U_{j_*}^\kappa(x),\\
 & \max\{P(x,t),U_{j_*}^\kappa(x)+\tilde{w}^\kappa(x)\}
 \le U_{j_*}^{\kappa+\epsilon}(x)+\tilde{w}^\kappa(x)+a\phi_1^\kappa(x)
 \le CU_{j_*}^\kappa(x),\\
 & \min\{P(x,t),U_{j_*}^\kappa(x)+\tilde{w}^\kappa(x)\}\ge U_{j_*}^\kappa(x).
\end{split}
\end{equation*}
These imply that 
\begin{equation}
\label{eq:6.7}
\begin{split}
0 & \le P(x,t)^{p-1}-\left[U_{j_*}^{\kappa}(x)+\tilde{w}^{\kappa}(x)\right]^{p-1}\\
 & \le C(\epsilon+a)(U_{j_*}^\kappa(x))^{p-1}
\le C(\epsilon+a)(u^\kappa(x))^{p-1}
\end{split}
\end{equation}
for $x\in{\bf R}^N$ and $t\in(0,1)$.

Let $\delta\in(0,1)$. 
By \eqref{eq:6.2}, \eqref{eq:6.5}, \eqref{eq:6.6} and \eqref{eq:6.7}, 
taking sufficiently small $a\in(0,1]$ and $\epsilon\in(0,1]$, 
we obtain 
\begin{equation*}
\begin{split}
H'(x,t) & \le ap(u^\kappa(x))^{p-1}\phi_1^\kappa(x)\\
 & \qquad
 +C\epsilon^{\min\{p-1,1\}}(u^\kappa(x))^{p-1}\phi_1^\kappa(x)+Ca(\epsilon+a)(u^\kappa(x))^{p-1}\phi_1^\kappa(x)\\
 & \le ap(1+\delta)(u^\kappa(x))^{p-1}\phi_1^\kappa(x)
\end{split}
\end{equation*}
for $x\in{\bf R}^N$ and $t\in(0,1)$.
This together with \eqref{eq:6.1} implies that
\begin{equation}
\label{eq:6.8}
\overline{u}-G*\overline{u}^p-(\kappa+\epsilon)G*\mu
\ge -ap(1+\delta)G*[(u^\kappa)^{p-1}\phi_1^\kappa]+a\phi_1^\kappa
\quad\mbox{in}\quad{\bf R}^N.
\end{equation}
On the other hand, it follows from $(\mbox{E}_\kappa)$ that
\begin{equation}
\label{eq:6.9}
\phi_1^\kappa=\lambda_1^\kappa p[G*(u^\kappa)^{p-1}\phi_1^\kappa]\quad\mbox{in}\quad{\bf R}^N.
\end{equation}
Combing \eqref{eq:6.8} with \eqref{eq:6.9}, we see that 
$$
\overline{u}-G*\overline{u}^p-(\kappa+\epsilon)G*\mu\ge a\left(1-\frac{1+\delta}{\lambda_1^\kappa}\right)\phi_1^\kappa
\quad\mbox{in}\quad{\bf R}^N.
$$
Since $\lambda_1^\kappa>1$ and $\phi_1^\kappa>0$ in ${\bf R}^N$, 
taking a sufficiently small $\delta\in(0,1)$ if necessary, 
we obtain 
$$
\overline{u}-G*\overline{u}^p-(\kappa+\epsilon)G*\mu>0
\quad\mbox{in}\quad{\bf R}^N.
$$
This means that $\overline{u}$ is a $(C_0+L^q_{\rm c})$-supersolution of \eqref{eq:1.1} 
with $\kappa$ replaced by $\kappa+\epsilon$. 
Then, by Lemma~\ref{Lemma:4.3} 
there exists a $(C_0+L^q_{\rm c})$-solution of \eqref{eq:1.1} with $\kappa$ replaced by $\kappa+\epsilon$. 
Therefore we see that $\kappa<\kappa^*$ and the proof is complete.
$\Box$
\vspace{5pt}

Now we are ready to prove Theorem~\ref{Theorem:1.2}.
\vspace{5pt}\newline
{\bf Proof of Theorem~\ref{Theorem:1.2}.}
By Proposition~\ref{Proposition:5.1} 
it suffices to prove the uniqueness of $(C_0+L^q_{\rm c})$-solution of \eqref{eq:1.1} 
with $\kappa=\kappa^*$. 

Let $u^{\kappa^*}$ be a minimal $(C_0+L^q_{\rm c})$-solution of \eqref{eq:1.1} with $\kappa=\kappa^*$ 
and set $w^{\kappa^*}=u^{\kappa^*}-U_{j_*}^{\kappa^*}$. 
Let $\tilde{u}$ be a $(C_{0}+L^{p}_{\rm c})$ solution of \eqref{eq:1.1} with $\kappa=\kappa^*$ and set $\tilde{w}=\tilde{u}-U_{j_*}^{\kappa^*}$. 
Then $w^{\kappa^*}$ and $\tilde{w}$ are solutions of \eqref{eq:4.3} with $\kappa=\kappa^*$
and the function $z:=\tilde{w}-w^{\kappa^*}$ satisfies
\begin{equation}
\label{eq:6.10}
z\ge 0\quad\mbox{in}\quad{\bf R}^N,
\quad
-\Delta z+z=(\tilde{w}+U_{j_*}^{\kappa^*})^p-(w^{\kappa^*}+U_{j_*}^{\kappa^*})^p
\quad\mbox{in}\quad{\bf R}^N.
\end{equation}

Let $\phi_1^{\kappa^*}$ be the eigenfunction of problem~$(\mbox{E}_{\kappa^{*}})$ corresponding to $\lambda_1^{\kappa^*}$ 
such that $\phi_1^{\kappa^*}>0$ in ${\bf R}^N$. 
By Lemmas~\ref{Lemma:6.1} and \ref{Lemma:6.2} we see that $\lambda_1^{\kappa^*}=1$. 
Then 
\begin{equation}
\label{eq:6.11}
-\Delta\phi_1^{\kappa^*}+\phi_1^{\kappa^*}=p(w^{\kappa^*}+U_{j_*}^{\kappa^*})^{p-1}\phi_1^{\kappa^*}\quad\mbox{in}\quad{\bf R}^N.
\end{equation}
Multiplying \eqref{eq:6.11} by $z$ and integrating it on ${\bf R}^N$, 
we obtain 
$$
\int_{{\bf R}^N}[\nabla z\nabla\phi_1^{\kappa^*}+z\phi_1^{\kappa^*}]\,dx
=p\int_{{\bf R}^N}(w^{\kappa^*}+U_{j_*}^{\kappa^*})^{p-1}\phi_1^{\kappa^*}(\tilde{w}-w^{\kappa^*})\,dx.
$$
On the other hand, 
by \eqref{eq:6.10} we see that 
$$
\int_{{\bf R}^N}[\nabla z\nabla\phi_1^{\kappa^*}+z\phi_1^{\kappa^*}]\,dx
=\int_{{\bf R}^N}
[(\tilde{w}+U_{j_*}^{\kappa^*})^p-(w^{\kappa^*}+U_{j_*}^{\kappa^*})^p]\phi_1^{\kappa^*}\,dx.
$$
These imply that 
$$
\int_{{\bf R}^N}[(\tilde{w}+U_{j_*}^{\kappa^*})^p-(w^{\kappa^*}+U_{j_*}^{\kappa^*})^p
-p(w^{\kappa^*}+U_{j_*}^{\kappa^*})^{p-1}(\tilde{w}-w^{\kappa^*})]\phi_1^{\kappa^*}\,dx=0.
$$
Since $t^p>s^p+ps^{p-1}(t-s)$ for $t>s\ge 0$, we obtain  
$\tilde{w}=w^{\kappa^*}$ in ${\bf R}^N$. 
Therefore we deduce that $u^{\kappa^*}$ is the unique $(C_0+L^q_{\rm c})$-solution of \eqref{eq:1.1} with $\kappa=\kappa^*$. 
Thus Theorem~\ref{Theorem:1.2} follows.
$\Box$ 
\vspace{5pt}

\noindent
{\bf Acknowledgements.} 
The first author was partially supported 
by the Grant-in-Aid for Scientific Research (A)(No.~15H02058)
from Japan Society for the Promotion of Science. 
%%%%%%%%%%%%%%%%%%%%%%%%%%%%%%%%%%%%%%%%%%%%%%%%
%%%%%%%%%%%%    references    %%%%%%%%%%%%%%%%%%
%%%%%%%%%%%%%%%%%%%%%%%%%%%%%%%%%%%%%%%%%%%%%%%%

\end{document}